\documentclass[11pt]{article}
\usepackage[margin=1in]{geometry}
\makeatletter
\renewcommand\@biblabel[1]{}
\renewenvironment{thebibliography}[1]
     {\section*{\refname}%
      \@mkboth{\MakeUppercase\refname}{\MakeUppercase\refname}%
      \list{}%
           {\leftmargin0pt
            \@openbib@code
            \usecounter{enumiv}}%
      \sloppy
      \clubpenalty4000
      \@clubpenalty \clubpenalty
      \widowpenalty4000%
      \sfcode`\.\@m}
     {\def\@noitemerr
       {\@latex@warning{Empty `thebibliography' environment}}%
      \endlist}
\makeatother

\usepackage{color}
\definecolor{AliceBlue}{rgb}{0.94,0.97,1.00}
\definecolor{AntiqueWhite1}{rgb}{1.00,0.94,0.86}
\definecolor{AntiqueWhite2}{rgb}{0.93,0.87,0.80}
\definecolor{AntiqueWhite3}{rgb}{0.80,0.75,0.69}
\definecolor{AntiqueWhite4}{rgb}{0.55,0.51,0.47}
\definecolor{AntiqueWhite}{rgb}{0.98,0.92,0.84}
\definecolor{BlanchedAlmond}{rgb}{1.00,0.92,0.80}
\definecolor{BlueViolet}{rgb}{0.54,0.17,0.89}
\definecolor{CadetBlue1}{rgb}{0.60,0.96,1.00}
\definecolor{CadetBlue2}{rgb}{0.56,0.90,0.93}
\definecolor{CadetBlue3}{rgb}{0.48,0.77,0.80}
\definecolor{CadetBlue4}{rgb}{0.33,0.53,0.55}
\definecolor{CadetBlue}{rgb}{0.37,0.62,0.63}
\definecolor{CornflowerBlue}{rgb}{0.39,0.58,0.93}
\definecolor{DarkBlue}{rgb}{0.00,0.00,0.55}
\definecolor{DarkCyan}{rgb}{0.00,0.55,0.55}
\definecolor{DarkGoldenrod1}{rgb}{1.00,0.73,0.06}
\definecolor{DarkGoldenrod2}{rgb}{0.93,0.68,0.05}
\definecolor{DarkGoldenrod3}{rgb}{0.80,0.58,0.05}
\definecolor{DarkGoldenrod4}{rgb}{0.55,0.40,0.03}
\definecolor{DarkGoldenrod}{rgb}{0.72,0.53,0.04}
\definecolor{DarkGray}{rgb}{0.66,0.66,0.66}
\definecolor{DarkGreen}{rgb}{0.00,0.39,0.00}
\definecolor{DarkGrey}{rgb}{0.66,0.66,0.66}
\definecolor{DarkKhaki}{rgb}{0.74,0.72,0.42}
\definecolor{DarkMagenta}{rgb}{0.55,0.00,0.55}
\definecolor{DarkOliveGreen1}{rgb}{0.79,1.00,0.44}
\definecolor{DarkOliveGreen2}{rgb}{0.74,0.93,0.41}
\definecolor{DarkOliveGreen3}{rgb}{0.64,0.80,0.35}
\definecolor{DarkOliveGreen4}{rgb}{0.43,0.55,0.24}
\definecolor{DarkOliveGreen}{rgb}{0.33,0.42,0.18}
\definecolor{DarkOrange1}{rgb}{1.00,0.50,0.00}
\definecolor{DarkOrange2}{rgb}{0.93,0.46,0.00}
\definecolor{DarkOrange3}{rgb}{0.80,0.40,0.00}
\definecolor{DarkOrange4}{rgb}{0.55,0.27,0.00}
\definecolor{DarkOrange}{rgb}{1.00,0.55,0.00}
\definecolor{DarkOrchid1}{rgb}{0.75,0.24,1.00}
\definecolor{DarkOrchid2}{rgb}{0.70,0.23,0.93}
\definecolor{DarkOrchid3}{rgb}{0.60,0.20,0.80}
\definecolor{DarkOrchid4}{rgb}{0.41,0.13,0.55}
\definecolor{DarkOrchid}{rgb}{0.60,0.20,0.80}
\definecolor{DarkRed}{rgb}{0.55,0.00,0.00}
\definecolor{DarkSalmon}{rgb}{0.91,0.59,0.48}
\definecolor{DarkSeaGreen1}{rgb}{0.76,1.00,0.76}
\definecolor{DarkSeaGreen2}{rgb}{0.71,0.93,0.71}
\definecolor{DarkSeaGreen3}{rgb}{0.61,0.80,0.61}
\definecolor{DarkSeaGreen4}{rgb}{0.41,0.55,0.41}
\definecolor{DarkSeaGreen}{rgb}{0.56,0.74,0.56}
\definecolor{DarkSlateBlue}{rgb}{0.28,0.24,0.55}
\definecolor{DarkSlateGray1}{rgb}{0.59,1.00,1.00}
\definecolor{DarkSlateGray2}{rgb}{0.55,0.93,0.93}
\definecolor{DarkSlateGray3}{rgb}{0.47,0.80,0.80}
\definecolor{DarkSlateGray4}{rgb}{0.32,0.55,0.55}
\definecolor{DarkSlateGray}{rgb}{0.18,0.31,0.31}
\definecolor{DarkSlateGrey}{rgb}{0.18,0.31,0.31}
\definecolor{DarkTurquoise}{rgb}{0.00,0.81,0.82}
\definecolor{DarkViolet}{rgb}{0.58,0.00,0.83}
\definecolor{DeepPink1}{rgb}{1.00,0.08,0.58}
\definecolor{DeepPink2}{rgb}{0.93,0.07,0.54}
\definecolor{DeepPink3}{rgb}{0.80,0.06,0.46}
\definecolor{DeepPink4}{rgb}{0.55,0.04,0.31}
\definecolor{DeepPink}{rgb}{1.00,0.08,0.58}
\definecolor{DeepSkyBlue1}{rgb}{0.00,0.75,1.00}
\definecolor{DeepSkyBlue2}{rgb}{0.00,0.70,0.93}
\definecolor{DeepSkyBlue3}{rgb}{0.00,0.60,0.80}
\definecolor{DeepSkyBlue4}{rgb}{0.00,0.41,0.55}
\definecolor{DeepSkyBlue}{rgb}{0.00,0.75,1.00}
\definecolor{DimGray}{rgb}{0.41,0.41,0.41}
\definecolor{DimGrey}{rgb}{0.41,0.41,0.41}
\definecolor{DodgerBlue1}{rgb}{0.12,0.56,1.00}
\definecolor{DodgerBlue2}{rgb}{0.11,0.53,0.93}
\definecolor{DodgerBlue3}{rgb}{0.09,0.45,0.80}
\definecolor{DodgerBlue4}{rgb}{0.06,0.31,0.55}
\definecolor{DodgerBlue}{rgb}{0.12,0.56,1.00}
\definecolor{FloralWhite}{rgb}{1.00,0.98,0.94}
\definecolor{ForestGreen}{rgb}{0.13,0.55,0.13}
\definecolor{GhostWhite}{rgb}{0.97,0.97,1.00}
\definecolor{GreenYellow}{rgb}{0.68,1.00,0.18}
\definecolor{HotPink1}{rgb}{1.00,0.43,0.71}
\definecolor{HotPink2}{rgb}{0.93,0.42,0.65}
\definecolor{HotPink3}{rgb}{0.80,0.38,0.56}
\definecolor{HotPink4}{rgb}{0.55,0.23,0.38}
\definecolor{HotPink}{rgb}{1.00,0.41,0.71}
\definecolor{IndianRed1}{rgb}{1.00,0.42,0.42}
\definecolor{IndianRed2}{rgb}{0.93,0.39,0.39}
\definecolor{IndianRed3}{rgb}{0.80,0.33,0.33}
\definecolor{IndianRed4}{rgb}{0.55,0.23,0.23}
\definecolor{IndianRed}{rgb}{0.80,0.36,0.36}
\definecolor{LavenderBlush1}{rgb}{1.00,0.94,0.96}
\definecolor{LavenderBlush2}{rgb}{0.93,0.88,0.90}
\definecolor{LavenderBlush3}{rgb}{0.80,0.76,0.77}
\definecolor{LavenderBlush4}{rgb}{0.55,0.51,0.53}
\definecolor{LavenderBlush}{rgb}{1.00,0.94,0.96}
\definecolor{LawnGreen}{rgb}{0.49,0.99,0.00}
\definecolor{LemonChiffon1}{rgb}{1.00,0.98,0.80}
\definecolor{LemonChiffon2}{rgb}{0.93,0.91,0.75}
\definecolor{LemonChiffon3}{rgb}{0.80,0.79,0.65}
\definecolor{LemonChiffon4}{rgb}{0.55,0.54,0.44}
\definecolor{LemonChiffon}{rgb}{1.00,0.98,0.80}
\definecolor{LightBlue1}{rgb}{0.75,0.94,1.00}
\definecolor{LightBlue2}{rgb}{0.70,0.87,0.93}
\definecolor{LightBlue3}{rgb}{0.60,0.75,0.80}
\definecolor{LightBlue4}{rgb}{0.41,0.51,0.55}
\definecolor{LightBlue}{rgb}{0.68,0.85,0.90}
\definecolor{LightCoral}{rgb}{0.94,0.50,0.50}
\definecolor{LightCyan1}{rgb}{0.88,1.00,1.00}
\definecolor{LightCyan2}{rgb}{0.82,0.93,0.93}
\definecolor{LightCyan3}{rgb}{0.71,0.80,0.80}
\definecolor{LightCyan4}{rgb}{0.48,0.55,0.55}
\definecolor{LightCyan}{rgb}{0.88,1.00,1.00}
\definecolor{LightGoldenrod1}{rgb}{1.00,0.93,0.55}
\definecolor{LightGoldenrod2}{rgb}{0.93,0.86,0.51}
\definecolor{LightGoldenrod3}{rgb}{0.80,0.75,0.44}
\definecolor{LightGoldenrod4}{rgb}{0.55,0.51,0.30}
\definecolor{LightGoldenrodYellow}{rgb}{0.98,0.98,0.82}
\definecolor{LightGoldenrod}{rgb}{0.93,0.87,0.51}
\definecolor{LightGray}{rgb}{0.83,0.83,0.83}
\definecolor{LightGreen}{rgb}{0.56,0.93,0.56}
\definecolor{LightGrey}{rgb}{0.83,0.83,0.83}
\definecolor{LightPink1}{rgb}{1.00,0.68,0.73}
\definecolor{LightPink2}{rgb}{0.93,0.64,0.68}
\definecolor{LightPink3}{rgb}{0.80,0.55,0.58}
\definecolor{LightPink4}{rgb}{0.55,0.37,0.40}
\definecolor{LightPink}{rgb}{1.00,0.71,0.76}
\definecolor{LightSalmon1}{rgb}{1.00,0.63,0.48}
\definecolor{LightSalmon2}{rgb}{0.93,0.58,0.45}
\definecolor{LightSalmon3}{rgb}{0.80,0.51,0.38}
\definecolor{LightSalmon4}{rgb}{0.55,0.34,0.26}
\definecolor{LightSalmon}{rgb}{1.00,0.63,0.48}
\definecolor{LightSeaGreen}{rgb}{0.13,0.70,0.67}
\definecolor{LightSkyBlue1}{rgb}{0.69,0.89,1.00}
\definecolor{LightSkyBlue2}{rgb}{0.64,0.83,0.93}
\definecolor{LightSkyBlue3}{rgb}{0.55,0.71,0.80}
\definecolor{LightSkyBlue4}{rgb}{0.38,0.48,0.55}
\definecolor{LightSkyBlue}{rgb}{0.53,0.81,0.98}
\definecolor{LightSlateBlue}{rgb}{0.52,0.44,1.00}
\definecolor{LightSlateGray}{rgb}{0.47,0.53,0.60}
\definecolor{LightSlateGrey}{rgb}{0.47,0.53,0.60}
\definecolor{LightSteelBlue1}{rgb}{0.79,0.88,1.00}
\definecolor{LightSteelBlue2}{rgb}{0.74,0.82,0.93}
\definecolor{LightSteelBlue3}{rgb}{0.64,0.71,0.80}
\definecolor{LightSteelBlue4}{rgb}{0.43,0.48,0.55}
\definecolor{LightSteelBlue}{rgb}{0.69,0.77,0.87}
\definecolor{LightYellow1}{rgb}{1.00,1.00,0.88}
\definecolor{LightYellow2}{rgb}{0.93,0.93,0.82}
\definecolor{LightYellow3}{rgb}{0.80,0.80,0.71}
\definecolor{LightYellow4}{rgb}{0.55,0.55,0.48}
\definecolor{LightYellow}{rgb}{1.00,1.00,0.88}
\definecolor{LimeGreen}{rgb}{0.20,0.80,0.20}
\definecolor{MediumAquamarine}{rgb}{0.40,0.80,0.67}
\definecolor{MediumBlue}{rgb}{0.00,0.00,0.80}
\definecolor{MediumOrchid1}{rgb}{0.88,0.40,1.00}
\definecolor{MediumOrchid2}{rgb}{0.82,0.37,0.93}
\definecolor{MediumOrchid3}{rgb}{0.71,0.32,0.80}
\definecolor{MediumOrchid4}{rgb}{0.48,0.22,0.55}
\definecolor{MediumOrchid}{rgb}{0.73,0.33,0.83}
\definecolor{MediumPurple1}{rgb}{0.67,0.51,1.00}
\definecolor{MediumPurple2}{rgb}{0.62,0.47,0.93}
\definecolor{MediumPurple3}{rgb}{0.54,0.41,0.80}
\definecolor{MediumPurple4}{rgb}{0.36,0.28,0.55}
\definecolor{MediumPurple}{rgb}{0.58,0.44,0.86}
\definecolor{MediumSeaGreen}{rgb}{0.24,0.70,0.44}
\definecolor{MediumSlateBlue}{rgb}{0.48,0.41,0.93}
\definecolor{MediumSpringGreen}{rgb}{0.00,0.98,0.60}
\definecolor{MediumTurquoise}{rgb}{0.28,0.82,0.80}
\definecolor{MediumVioletRed}{rgb}{0.78,0.08,0.52}
\definecolor{MidnightBlue}{rgb}{0.10,0.10,0.44}
\definecolor{MintCream}{rgb}{0.96,1.00,0.98}
\definecolor{MistyRose1}{rgb}{1.00,0.89,0.88}
\definecolor{MistyRose2}{rgb}{0.93,0.84,0.82}
\definecolor{MistyRose3}{rgb}{0.80,0.72,0.71}
\definecolor{MistyRose4}{rgb}{0.55,0.49,0.48}
\definecolor{MistyRose}{rgb}{1.00,0.89,0.88}
\definecolor{NavajoWhite1}{rgb}{1.00,0.87,0.68}
\definecolor{NavajoWhite2}{rgb}{0.93,0.81,0.63}
\definecolor{NavajoWhite3}{rgb}{0.80,0.70,0.55}
\definecolor{NavajoWhite4}{rgb}{0.55,0.47,0.37}
\definecolor{NavajoWhite}{rgb}{1.00,0.87,0.68}
\definecolor{NavyBlue}{rgb}{0.00,0.00,0.50}
\definecolor{OldLace}{rgb}{0.99,0.96,0.90}
\definecolor{OliveDrab1}{rgb}{0.75,1.00,0.24}
\definecolor{OliveDrab2}{rgb}{0.70,0.93,0.23}
\definecolor{OliveDrab3}{rgb}{0.60,0.80,0.20}
\definecolor{OliveDrab4}{rgb}{0.41,0.55,0.13}
\definecolor{OliveDrab}{rgb}{0.42,0.56,0.14}
\definecolor{OrangeRed1}{rgb}{1.00,0.27,0.00}
\definecolor{OrangeRed2}{rgb}{0.93,0.25,0.00}
\definecolor{OrangeRed3}{rgb}{0.80,0.22,0.00}
\definecolor{OrangeRed4}{rgb}{0.55,0.15,0.00}
\definecolor{OrangeRed}{rgb}{1.00,0.27,0.00}
\definecolor{PaleGoldenrod}{rgb}{0.93,0.91,0.67}
\definecolor{PaleGreen1}{rgb}{0.60,1.00,0.60}
\definecolor{PaleGreen2}{rgb}{0.56,0.93,0.56}
\definecolor{PaleGreen3}{rgb}{0.49,0.80,0.49}
\definecolor{PaleGreen4}{rgb}{0.33,0.55,0.33}
\definecolor{PaleGreen}{rgb}{0.60,0.98,0.60}
\definecolor{PaleTurquoise1}{rgb}{0.73,1.00,1.00}
\definecolor{PaleTurquoise2}{rgb}{0.68,0.93,0.93}
\definecolor{PaleTurquoise3}{rgb}{0.59,0.80,0.80}
\definecolor{PaleTurquoise4}{rgb}{0.40,0.55,0.55}
\definecolor{PaleTurquoise}{rgb}{0.69,0.93,0.93}
\definecolor{PaleVioletRed1}{rgb}{1.00,0.51,0.67}
\definecolor{PaleVioletRed2}{rgb}{0.93,0.47,0.62}
\definecolor{PaleVioletRed3}{rgb}{0.80,0.41,0.54}
\definecolor{PaleVioletRed4}{rgb}{0.55,0.28,0.36}
\definecolor{PaleVioletRed}{rgb}{0.86,0.44,0.58}
\definecolor{PapayaWhip}{rgb}{1.00,0.94,0.84}
\definecolor{PeachPuff1}{rgb}{1.00,0.85,0.73}
\definecolor{PeachPuff2}{rgb}{0.93,0.80,0.68}
\definecolor{PeachPuff3}{rgb}{0.80,0.69,0.58}
\definecolor{PeachPuff4}{rgb}{0.55,0.47,0.40}
\definecolor{PeachPuff}{rgb}{1.00,0.85,0.73}
\definecolor{PowderBlue}{rgb}{0.69,0.88,0.90}
\definecolor{RosyBrown1}{rgb}{1.00,0.76,0.76}
\definecolor{RosyBrown2}{rgb}{0.93,0.71,0.71}
\definecolor{RosyBrown3}{rgb}{0.80,0.61,0.61}
\definecolor{RosyBrown4}{rgb}{0.55,0.41,0.41}
\definecolor{RosyBrown}{rgb}{0.74,0.56,0.56}
\definecolor{RoyalBlue1}{rgb}{0.28,0.46,1.00}
\definecolor{RoyalBlue2}{rgb}{0.26,0.43,0.93}
\definecolor{RoyalBlue3}{rgb}{0.23,0.37,0.80}
\definecolor{RoyalBlue4}{rgb}{0.15,0.25,0.55}
\definecolor{RoyalBlue}{rgb}{0.25,0.41,0.88}
\definecolor{SaddleBrown}{rgb}{0.55,0.27,0.07}
\definecolor{SandyBrown}{rgb}{0.96,0.64,0.38}
\definecolor{SeaGreen1}{rgb}{0.33,1.00,0.62}
\definecolor{SeaGreen2}{rgb}{0.31,0.93,0.58}
\definecolor{SeaGreen3}{rgb}{0.26,0.80,0.50}
\definecolor{SeaGreen4}{rgb}{0.18,0.55,0.34}
\definecolor{SeaGreen}{rgb}{0.18,0.55,0.34}
\definecolor{SkyBlue1}{rgb}{0.53,0.81,1.00}
\definecolor{SkyBlue2}{rgb}{0.49,0.75,0.93}
\definecolor{SkyBlue3}{rgb}{0.42,0.65,0.80}
\definecolor{SkyBlue4}{rgb}{0.29,0.44,0.55}
\definecolor{SkyBlue}{rgb}{0.53,0.81,0.92}
\definecolor{SlateBlue1}{rgb}{0.51,0.44,1.00}
\definecolor{SlateBlue2}{rgb}{0.48,0.40,0.93}
\definecolor{SlateBlue3}{rgb}{0.41,0.35,0.80}
\definecolor{SlateBlue4}{rgb}{0.28,0.24,0.55}
\definecolor{SlateBlue}{rgb}{0.42,0.35,0.80}
\definecolor{SlateGray1}{rgb}{0.78,0.89,1.00}
\definecolor{SlateGray2}{rgb}{0.73,0.83,0.93}
\definecolor{SlateGray3}{rgb}{0.62,0.71,0.80}
\definecolor{SlateGray4}{rgb}{0.42,0.48,0.55}
\definecolor{SlateGray}{rgb}{0.44,0.50,0.56}
\definecolor{SlateGrey}{rgb}{0.44,0.50,0.56}
\definecolor{SpringGreen1}{rgb}{0.00,1.00,0.50}
\definecolor{SpringGreen2}{rgb}{0.00,0.93,0.46}
\definecolor{SpringGreen3}{rgb}{0.00,0.80,0.40}
\definecolor{SpringGreen4}{rgb}{0.00,0.55,0.27}
\definecolor{SpringGreen}{rgb}{0.00,1.00,0.50}
\definecolor{SteelBlue1}{rgb}{0.39,0.72,1.00}
\definecolor{SteelBlue2}{rgb}{0.36,0.67,0.93}
\definecolor{SteelBlue3}{rgb}{0.31,0.58,0.80}
\definecolor{SteelBlue4}{rgb}{0.21,0.39,0.55}
\definecolor{SteelBlue}{rgb}{0.27,0.51,0.71}
\definecolor{VioletRed1}{rgb}{1.00,0.24,0.59}
\definecolor{VioletRed2}{rgb}{0.93,0.23,0.55}
\definecolor{VioletRed3}{rgb}{0.80,0.20,0.47}
\definecolor{VioletRed4}{rgb}{0.55,0.13,0.32}
\definecolor{VioletRed}{rgb}{0.82,0.13,0.56}
\definecolor{WhiteSmoke}{rgb}{0.96,0.96,0.96}
\definecolor{YellowGreen}{rgb}{0.60,0.80,0.20}
\definecolor{aliceblue}{rgb}{0.94,0.97,1.00}
\definecolor{antiquewhite}{rgb}{0.98,0.92,0.84}
\definecolor{aquamarine1}{rgb}{0.50,1.00,0.83}
\definecolor{aquamarine2}{rgb}{0.46,0.93,0.78}
\definecolor{aquamarine3}{rgb}{0.40,0.80,0.67}
\definecolor{aquamarine4}{rgb}{0.27,0.55,0.45}
\definecolor{aquamarine}{rgb}{0.50,1.00,0.83}
\definecolor{azure1}{rgb}{0.94,1.00,1.00}
\definecolor{azure2}{rgb}{0.88,0.93,0.93}
\definecolor{azure3}{rgb}{0.76,0.80,0.80}
\definecolor{azure4}{rgb}{0.51,0.55,0.55}
\definecolor{azure}{rgb}{0.94,1.00,1.00}
\definecolor{beige}{rgb}{0.96,0.96,0.86}
\definecolor{bisque1}{rgb}{1.00,0.89,0.77}
\definecolor{bisque2}{rgb}{0.93,0.84,0.72}
\definecolor{bisque3}{rgb}{0.80,0.72,0.62}
\definecolor{bisque4}{rgb}{0.55,0.49,0.42}
\definecolor{bisque}{rgb}{1.00,0.89,0.77}
\definecolor{black}{rgb}{0.00,0.00,0.00}
\definecolor{blanchedalmond}{rgb}{1.00,0.92,0.80}
\definecolor{blue1}{rgb}{0.00,0.00,1.00}
\definecolor{blue2}{rgb}{0.00,0.00,0.93}
\definecolor{blue3}{rgb}{0.00,0.00,0.80}
\definecolor{blue4}{rgb}{0.00,0.00,0.55}
\definecolor{blueviolet}{rgb}{0.54,0.17,0.89}
\definecolor{blue}{rgb}{0.00,0.00,1.00}
\definecolor{brown1}{rgb}{1.00,0.25,0.25}
\definecolor{brown2}{rgb}{0.93,0.23,0.23}
\definecolor{brown3}{rgb}{0.80,0.20,0.20}
\definecolor{brown4}{rgb}{0.55,0.14,0.14}
\definecolor{brown}{rgb}{0.65,0.16,0.16}
\definecolor{burlywood1}{rgb}{1.00,0.83,0.61}
\definecolor{burlywood2}{rgb}{0.93,0.77,0.57}
\definecolor{burlywood3}{rgb}{0.80,0.67,0.49}
\definecolor{burlywood4}{rgb}{0.55,0.45,0.33}
\definecolor{burlywood}{rgb}{0.87,0.72,0.53}
\definecolor{cadetblue}{rgb}{0.37,0.62,0.63}
\definecolor{chartreuse1}{rgb}{0.50,1.00,0.00}
\definecolor{chartreuse2}{rgb}{0.46,0.93,0.00}
\definecolor{chartreuse3}{rgb}{0.40,0.80,0.00}
\definecolor{chartreuse4}{rgb}{0.27,0.55,0.00}
\definecolor{chartreuse}{rgb}{0.50,1.00,0.00}
\definecolor{chocolate1}{rgb}{1.00,0.50,0.14}
\definecolor{chocolate2}{rgb}{0.93,0.46,0.13}
\definecolor{chocolate3}{rgb}{0.80,0.40,0.11}
\definecolor{chocolate4}{rgb}{0.55,0.27,0.07}
\definecolor{chocolate}{rgb}{0.82,0.41,0.12}
\definecolor{coral1}{rgb}{1.00,0.45,0.34}
\definecolor{coral2}{rgb}{0.93,0.42,0.31}
\definecolor{coral3}{rgb}{0.80,0.36,0.27}
\definecolor{coral4}{rgb}{0.55,0.24,0.18}
\definecolor{coral}{rgb}{1.00,0.50,0.31}
\definecolor{cornflowerblue}{rgb}{0.39,0.58,0.93}
\definecolor{cornsilk1}{rgb}{1.00,0.97,0.86}
\definecolor{cornsilk2}{rgb}{0.93,0.91,0.80}
\definecolor{cornsilk3}{rgb}{0.80,0.78,0.69}
\definecolor{cornsilk4}{rgb}{0.55,0.53,0.47}
\definecolor{cornsilk}{rgb}{1.00,0.97,0.86}
\definecolor{cyan1}{rgb}{0.00,1.00,1.00}
\definecolor{cyan2}{rgb}{0.00,0.93,0.93}
\definecolor{cyan3}{rgb}{0.00,0.80,0.80}
\definecolor{cyan4}{rgb}{0.00,0.55,0.55}
\definecolor{cyan}{rgb}{0.00,1.00,1.00}
\definecolor{darkblue}{rgb}{0.00,0.00,0.55}
\definecolor{darkcyan}{rgb}{0.00,0.55,0.55}
\definecolor{darkgoldenrod}{rgb}{0.72,0.53,0.04}
\definecolor{darkgray}{rgb}{0.66,0.66,0.66}
\definecolor{darkgreen}{rgb}{0.00,0.39,0.00}
\definecolor{darkgrey}{rgb}{0.66,0.66,0.66}
\definecolor{darkkhaki}{rgb}{0.74,0.72,0.42}
\definecolor{darkmagenta}{rgb}{0.55,0.00,0.55}
\definecolor{darkolive}{rgb}{0.33,0.42,0.18}
\definecolor{darkorange}{rgb}{1.00,0.55,0.00}
\definecolor{darkorchid}{rgb}{0.60,0.20,0.80}
\definecolor{darkred}{rgb}{0.55,0.00,0.00}
\definecolor{darksalmon}{rgb}{0.91,0.59,0.48}
\definecolor{darksea}{rgb}{0.56,0.74,0.56}
\definecolor{darkslate}{rgb}{0.18,0.31,0.31}
\definecolor{darkslate}{rgb}{0.18,0.31,0.31}
\definecolor{darkslate}{rgb}{0.28,0.24,0.55}
\definecolor{darkturquoise}{rgb}{0.00,0.81,0.82}
\definecolor{darkviolet}{rgb}{0.58,0.00,0.83}
\definecolor{deeppink}{rgb}{1.00,0.08,0.58}
\definecolor{deepsky}{rgb}{0.00,0.75,1.00}
\definecolor{dimgray}{rgb}{0.41,0.41,0.41}
\definecolor{dimgrey}{rgb}{0.41,0.41,0.41}
\definecolor{dodgerblue}{rgb}{0.12,0.56,1.00}
\definecolor{firebrick1}{rgb}{1.00,0.19,0.19}
\definecolor{firebrick2}{rgb}{0.93,0.17,0.17}
\definecolor{firebrick3}{rgb}{0.80,0.15,0.15}
\definecolor{firebrick4}{rgb}{0.55,0.10,0.10}
\definecolor{firebrick}{rgb}{0.70,0.13,0.13}
\definecolor{floralwhite}{rgb}{1.00,0.98,0.94}
\definecolor{forestgreen}{rgb}{0.13,0.55,0.13}
\definecolor{gainsboro}{rgb}{0.86,0.86,0.86}
\definecolor{ghostwhite}{rgb}{0.97,0.97,1.00}
\definecolor{gold1}{rgb}{1.00,0.84,0.00}
\definecolor{gold2}{rgb}{0.93,0.79,0.00}
\definecolor{gold3}{rgb}{0.80,0.68,0.00}
\definecolor{gold4}{rgb}{0.55,0.46,0.00}
\definecolor{goldenrod1}{rgb}{1.00,0.76,0.15}
\definecolor{goldenrod2}{rgb}{0.93,0.71,0.13}
\definecolor{goldenrod3}{rgb}{0.80,0.61,0.11}
\definecolor{goldenrod4}{rgb}{0.55,0.41,0.08}
\definecolor{goldenrod}{rgb}{0.85,0.65,0.13}
\definecolor{gold}{rgb}{1.00,0.84,0.00}
\definecolor{gray0}{rgb}{0.00,0.00,0.00}
\definecolor{gray100}{rgb}{1.00,1.00,1.00}
\definecolor{gray10}{rgb}{0.10,0.10,0.10}
\definecolor{gray11}{rgb}{0.11,0.11,0.11}
\definecolor{gray12}{rgb}{0.12,0.12,0.12}
\definecolor{gray13}{rgb}{0.13,0.13,0.13}
\definecolor{gray14}{rgb}{0.14,0.14,0.14}
\definecolor{gray15}{rgb}{0.15,0.15,0.15}
\definecolor{gray16}{rgb}{0.16,0.16,0.16}
\definecolor{gray17}{rgb}{0.17,0.17,0.17}
\definecolor{gray18}{rgb}{0.18,0.18,0.18}
\definecolor{gray19}{rgb}{0.19,0.19,0.19}
\definecolor{gray1}{rgb}{0.01,0.01,0.01}
\definecolor{gray20}{rgb}{0.20,0.20,0.20}
\definecolor{gray21}{rgb}{0.21,0.21,0.21}
\definecolor{gray22}{rgb}{0.22,0.22,0.22}
\definecolor{gray23}{rgb}{0.23,0.23,0.23}
\definecolor{gray24}{rgb}{0.24,0.24,0.24}
\definecolor{gray25}{rgb}{0.25,0.25,0.25}
\definecolor{gray26}{rgb}{0.26,0.26,0.26}
\definecolor{gray27}{rgb}{0.27,0.27,0.27}
\definecolor{gray28}{rgb}{0.28,0.28,0.28}
\definecolor{gray29}{rgb}{0.29,0.29,0.29}
\definecolor{gray2}{rgb}{0.02,0.02,0.02}
\definecolor{gray30}{rgb}{0.30,0.30,0.30}
\definecolor{gray31}{rgb}{0.31,0.31,0.31}
\definecolor{gray32}{rgb}{0.32,0.32,0.32}
\definecolor{gray33}{rgb}{0.33,0.33,0.33}
\definecolor{gray34}{rgb}{0.34,0.34,0.34}
\definecolor{gray35}{rgb}{0.35,0.35,0.35}
\definecolor{gray36}{rgb}{0.36,0.36,0.36}
\definecolor{gray37}{rgb}{0.37,0.37,0.37}
\definecolor{gray38}{rgb}{0.38,0.38,0.38}
\definecolor{gray39}{rgb}{0.39,0.39,0.39}
\definecolor{gray3}{rgb}{0.03,0.03,0.03}
\definecolor{gray40}{rgb}{0.40,0.40,0.40}
\definecolor{gray41}{rgb}{0.41,0.41,0.41}
\definecolor{gray42}{rgb}{0.42,0.42,0.42}
\definecolor{gray43}{rgb}{0.43,0.43,0.43}
\definecolor{gray44}{rgb}{0.44,0.44,0.44}
\definecolor{gray45}{rgb}{0.45,0.45,0.45}
\definecolor{gray46}{rgb}{0.46,0.46,0.46}
\definecolor{gray47}{rgb}{0.47,0.47,0.47}
\definecolor{gray48}{rgb}{0.48,0.48,0.48}
\definecolor{gray49}{rgb}{0.49,0.49,0.49}
\definecolor{gray4}{rgb}{0.04,0.04,0.04}
\definecolor{gray50}{rgb}{0.50,0.50,0.50}
\definecolor{gray51}{rgb}{0.51,0.51,0.51}
\definecolor{gray52}{rgb}{0.52,0.52,0.52}
\definecolor{gray53}{rgb}{0.53,0.53,0.53}
\definecolor{gray54}{rgb}{0.54,0.54,0.54}
\definecolor{gray55}{rgb}{0.55,0.55,0.55}
\definecolor{gray56}{rgb}{0.56,0.56,0.56}
\definecolor{gray57}{rgb}{0.57,0.57,0.57}
\definecolor{gray58}{rgb}{0.58,0.58,0.58}
\definecolor{gray59}{rgb}{0.59,0.59,0.59}
\definecolor{gray5}{rgb}{0.05,0.05,0.05}
\definecolor{gray60}{rgb}{0.60,0.60,0.60}
\definecolor{gray61}{rgb}{0.61,0.61,0.61}
\definecolor{gray62}{rgb}{0.62,0.62,0.62}
\definecolor{gray63}{rgb}{0.63,0.63,0.63}
\definecolor{gray64}{rgb}{0.64,0.64,0.64}
\definecolor{gray65}{rgb}{0.65,0.65,0.65}
\definecolor{gray66}{rgb}{0.66,0.66,0.66}
\definecolor{gray67}{rgb}{0.67,0.67,0.67}
\definecolor{gray68}{rgb}{0.68,0.68,0.68}
\definecolor{gray69}{rgb}{0.69,0.69,0.69}
\definecolor{gray6}{rgb}{0.06,0.06,0.06}
\definecolor{gray70}{rgb}{0.70,0.70,0.70}
\definecolor{gray71}{rgb}{0.71,0.71,0.71}
\definecolor{gray72}{rgb}{0.72,0.72,0.72}
\definecolor{gray73}{rgb}{0.73,0.73,0.73}
\definecolor{gray74}{rgb}{0.74,0.74,0.74}
\definecolor{gray75}{rgb}{0.75,0.75,0.75}
\definecolor{gray76}{rgb}{0.76,0.76,0.76}
\definecolor{gray77}{rgb}{0.77,0.77,0.77}
\definecolor{gray78}{rgb}{0.78,0.78,0.78}
\definecolor{gray79}{rgb}{0.79,0.79,0.79}
\definecolor{gray7}{rgb}{0.07,0.07,0.07}
\definecolor{gray80}{rgb}{0.80,0.80,0.80}
\definecolor{gray81}{rgb}{0.81,0.81,0.81}
\definecolor{gray82}{rgb}{0.82,0.82,0.82}
\definecolor{gray83}{rgb}{0.83,0.83,0.83}
\definecolor{gray84}{rgb}{0.84,0.84,0.84}
\definecolor{gray85}{rgb}{0.85,0.85,0.85}
\definecolor{gray86}{rgb}{0.86,0.86,0.86}
\definecolor{gray87}{rgb}{0.87,0.87,0.87}
\definecolor{gray88}{rgb}{0.88,0.88,0.88}
\definecolor{gray89}{rgb}{0.89,0.89,0.89}
\definecolor{gray8}{rgb}{0.08,0.08,0.08}
\definecolor{gray90}{rgb}{0.90,0.90,0.90}
\definecolor{gray91}{rgb}{0.91,0.91,0.91}
\definecolor{gray92}{rgb}{0.92,0.92,0.92}
\definecolor{gray93}{rgb}{0.93,0.93,0.93}
\definecolor{gray94}{rgb}{0.94,0.94,0.94}
\definecolor{gray95}{rgb}{0.95,0.95,0.95}
\definecolor{gray96}{rgb}{0.96,0.96,0.96}
\definecolor{gray97}{rgb}{0.97,0.97,0.97}
\definecolor{gray98}{rgb}{0.98,0.98,0.98}
\definecolor{gray99}{rgb}{0.99,0.99,0.99}
\definecolor{gray9}{rgb}{0.09,0.09,0.09}
\definecolor{gray}{rgb}{0.75,0.75,0.75}
\definecolor{green1}{rgb}{0.00,1.00,0.00}
\definecolor{green2}{rgb}{0.00,0.93,0.00}
\definecolor{green3}{rgb}{0.00,0.80,0.00}
\definecolor{green4}{rgb}{0.00,0.55,0.00}
\definecolor{greenyellow}{rgb}{0.68,1.00,0.18}
\definecolor{green}{rgb}{0.00,1.00,0.00}
\definecolor{grey0}{rgb}{0.00,0.00,0.00}
\definecolor{grey100}{rgb}{1.00,1.00,1.00}
\definecolor{grey10}{rgb}{0.10,0.10,0.10}
\definecolor{grey11}{rgb}{0.11,0.11,0.11}
\definecolor{grey12}{rgb}{0.12,0.12,0.12}
\definecolor{grey13}{rgb}{0.13,0.13,0.13}
\definecolor{grey14}{rgb}{0.14,0.14,0.14}
\definecolor{grey15}{rgb}{0.15,0.15,0.15}
\definecolor{grey16}{rgb}{0.16,0.16,0.16}
\definecolor{grey17}{rgb}{0.17,0.17,0.17}
\definecolor{grey18}{rgb}{0.18,0.18,0.18}
\definecolor{grey19}{rgb}{0.19,0.19,0.19}
\definecolor{grey1}{rgb}{0.01,0.01,0.01}
\definecolor{grey20}{rgb}{0.20,0.20,0.20}
\definecolor{grey21}{rgb}{0.21,0.21,0.21}
\definecolor{grey22}{rgb}{0.22,0.22,0.22}
\definecolor{grey23}{rgb}{0.23,0.23,0.23}
\definecolor{grey24}{rgb}{0.24,0.24,0.24}
\definecolor{grey25}{rgb}{0.25,0.25,0.25}
\definecolor{grey26}{rgb}{0.26,0.26,0.26}
\definecolor{grey27}{rgb}{0.27,0.27,0.27}
\definecolor{grey28}{rgb}{0.28,0.28,0.28}
\definecolor{grey29}{rgb}{0.29,0.29,0.29}
\definecolor{grey2}{rgb}{0.02,0.02,0.02}
\definecolor{grey30}{rgb}{0.30,0.30,0.30}
\definecolor{grey31}{rgb}{0.31,0.31,0.31}
\definecolor{grey32}{rgb}{0.32,0.32,0.32}
\definecolor{grey33}{rgb}{0.33,0.33,0.33}
\definecolor{grey34}{rgb}{0.34,0.34,0.34}
\definecolor{grey35}{rgb}{0.35,0.35,0.35}
\definecolor{grey36}{rgb}{0.36,0.36,0.36}
\definecolor{grey37}{rgb}{0.37,0.37,0.37}
\definecolor{grey38}{rgb}{0.38,0.38,0.38}
\definecolor{grey39}{rgb}{0.39,0.39,0.39}
\definecolor{grey3}{rgb}{0.03,0.03,0.03}
\definecolor{grey40}{rgb}{0.40,0.40,0.40}
\definecolor{grey41}{rgb}{0.41,0.41,0.41}
\definecolor{grey42}{rgb}{0.42,0.42,0.42}
\definecolor{grey43}{rgb}{0.43,0.43,0.43}
\definecolor{grey44}{rgb}{0.44,0.44,0.44}
\definecolor{grey45}{rgb}{0.45,0.45,0.45}
\definecolor{grey46}{rgb}{0.46,0.46,0.46}
\definecolor{grey47}{rgb}{0.47,0.47,0.47}
\definecolor{grey48}{rgb}{0.48,0.48,0.48}
\definecolor{grey49}{rgb}{0.49,0.49,0.49}
\definecolor{grey4}{rgb}{0.04,0.04,0.04}
\definecolor{grey50}{rgb}{0.50,0.50,0.50}
\definecolor{grey51}{rgb}{0.51,0.51,0.51}
\definecolor{grey52}{rgb}{0.52,0.52,0.52}
\definecolor{grey53}{rgb}{0.53,0.53,0.53}
\definecolor{grey54}{rgb}{0.54,0.54,0.54}
\definecolor{grey55}{rgb}{0.55,0.55,0.55}
\definecolor{grey56}{rgb}{0.56,0.56,0.56}
\definecolor{grey57}{rgb}{0.57,0.57,0.57}
\definecolor{grey58}{rgb}{0.58,0.58,0.58}
\definecolor{grey59}{rgb}{0.59,0.59,0.59}
\definecolor{grey5}{rgb}{0.05,0.05,0.05}
\definecolor{grey60}{rgb}{0.60,0.60,0.60}
\definecolor{grey61}{rgb}{0.61,0.61,0.61}
\definecolor{grey62}{rgb}{0.62,0.62,0.62}
\definecolor{grey63}{rgb}{0.63,0.63,0.63}
\definecolor{grey64}{rgb}{0.64,0.64,0.64}
\definecolor{grey65}{rgb}{0.65,0.65,0.65}
\definecolor{grey66}{rgb}{0.66,0.66,0.66}
\definecolor{grey67}{rgb}{0.67,0.67,0.67}
\definecolor{grey68}{rgb}{0.68,0.68,0.68}
\definecolor{grey69}{rgb}{0.69,0.69,0.69}
\definecolor{grey6}{rgb}{0.06,0.06,0.06}
\definecolor{grey70}{rgb}{0.70,0.70,0.70}
\definecolor{grey71}{rgb}{0.71,0.71,0.71}
\definecolor{grey72}{rgb}{0.72,0.72,0.72}
\definecolor{grey73}{rgb}{0.73,0.73,0.73}
\definecolor{grey74}{rgb}{0.74,0.74,0.74}
\definecolor{grey75}{rgb}{0.75,0.75,0.75}
\definecolor{grey76}{rgb}{0.76,0.76,0.76}
\definecolor{grey77}{rgb}{0.77,0.77,0.77}
\definecolor{grey78}{rgb}{0.78,0.78,0.78}
\definecolor{grey79}{rgb}{0.79,0.79,0.79}
\definecolor{grey7}{rgb}{0.07,0.07,0.07}
\definecolor{grey80}{rgb}{0.80,0.80,0.80}
\definecolor{grey81}{rgb}{0.81,0.81,0.81}
\definecolor{grey82}{rgb}{0.82,0.82,0.82}
\definecolor{grey83}{rgb}{0.83,0.83,0.83}
\definecolor{grey84}{rgb}{0.84,0.84,0.84}
\definecolor{grey85}{rgb}{0.85,0.85,0.85}
\definecolor{grey86}{rgb}{0.86,0.86,0.86}
\definecolor{grey87}{rgb}{0.87,0.87,0.87}
\definecolor{grey88}{rgb}{0.88,0.88,0.88}
\definecolor{grey89}{rgb}{0.89,0.89,0.89}
\definecolor{grey8}{rgb}{0.08,0.08,0.08}
\definecolor{grey90}{rgb}{0.90,0.90,0.90}
\definecolor{grey91}{rgb}{0.91,0.91,0.91}
\definecolor{grey92}{rgb}{0.92,0.92,0.92}
\definecolor{grey93}{rgb}{0.93,0.93,0.93}
\definecolor{grey94}{rgb}{0.94,0.94,0.94}
\definecolor{grey95}{rgb}{0.95,0.95,0.95}
\definecolor{grey96}{rgb}{0.96,0.96,0.96}
\definecolor{grey97}{rgb}{0.97,0.97,0.97}
\definecolor{grey98}{rgb}{0.98,0.98,0.98}
\definecolor{grey99}{rgb}{0.99,0.99,0.99}
\definecolor{grey9}{rgb}{0.09,0.09,0.09}
\definecolor{grey}{rgb}{0.75,0.75,0.75}
\definecolor{honeydew1}{rgb}{0.94,1.00,0.94}
\definecolor{honeydew2}{rgb}{0.88,0.93,0.88}
\definecolor{honeydew3}{rgb}{0.76,0.80,0.76}
\definecolor{honeydew4}{rgb}{0.51,0.55,0.51}
\definecolor{honeydew}{rgb}{0.94,1.00,0.94}
\definecolor{hotpink}{rgb}{1.00,0.41,0.71}
\definecolor{indianred}{rgb}{0.80,0.36,0.36}
\definecolor{ivory1}{rgb}{1.00,1.00,0.94}
\definecolor{ivory2}{rgb}{0.93,0.93,0.88}
\definecolor{ivory3}{rgb}{0.80,0.80,0.76}
\definecolor{ivory4}{rgb}{0.55,0.55,0.51}
\definecolor{ivory}{rgb}{1.00,1.00,0.94}
\definecolor{khaki1}{rgb}{1.00,0.96,0.56}
\definecolor{khaki2}{rgb}{0.93,0.90,0.52}
\definecolor{khaki3}{rgb}{0.80,0.78,0.45}
\definecolor{khaki4}{rgb}{0.55,0.53,0.31}
\definecolor{khaki}{rgb}{0.94,0.90,0.55}
\definecolor{lavenderblush}{rgb}{1.00,0.94,0.96}
\definecolor{lavender}{rgb}{0.90,0.90,0.98}
\definecolor{lawngreen}{rgb}{0.49,0.99,0.00}
\definecolor{lemonchiffon}{rgb}{1.00,0.98,0.80}
\definecolor{lightblue}{rgb}{0.68,0.85,0.90}
\definecolor{lightcoral}{rgb}{0.94,0.50,0.50}
\definecolor{lightcyan}{rgb}{0.88,1.00,1.00}
\definecolor{lightgoldenrod}{rgb}{0.93,0.87,0.51}
\definecolor{lightgoldenrod}{rgb}{0.98,0.98,0.82}
\definecolor{lightgray}{rgb}{0.83,0.83,0.83}
\definecolor{lightgreen}{rgb}{0.56,0.93,0.56}
\definecolor{lightgrey}{rgb}{0.83,0.83,0.83}
\definecolor{lightpink}{rgb}{1.00,0.71,0.76}
\definecolor{lightsalmon}{rgb}{1.00,0.63,0.48}
\definecolor{lightsea}{rgb}{0.13,0.70,0.67}
\definecolor{lightsky}{rgb}{0.53,0.81,0.98}
\definecolor{lightslate}{rgb}{0.47,0.53,0.60}
\definecolor{lightslate}{rgb}{0.47,0.53,0.60}
\definecolor{lightslate}{rgb}{0.52,0.44,1.00}
\definecolor{lightsteel}{rgb}{0.69,0.77,0.87}
\definecolor{lightyellow}{rgb}{1.00,1.00,0.88}
\definecolor{limegreen}{rgb}{0.20,0.80,0.20}
\definecolor{linen}{rgb}{0.98,0.94,0.90}
\definecolor{magenta1}{rgb}{1.00,0.00,1.00}
\definecolor{magenta2}{rgb}{0.93,0.00,0.93}
\definecolor{magenta3}{rgb}{0.80,0.00,0.80}
\definecolor{magenta4}{rgb}{0.55,0.00,0.55}
\definecolor{magenta}{rgb}{1.00,0.00,1.00}
\definecolor{maroon1}{rgb}{1.00,0.20,0.70}
\definecolor{maroon2}{rgb}{0.93,0.19,0.65}
\definecolor{maroon3}{rgb}{0.80,0.16,0.56}
\definecolor{maroon4}{rgb}{0.55,0.11,0.38}
\definecolor{maroon}{rgb}{0.69,0.19,0.38}
\definecolor{mediumaquamarine}{rgb}{0.40,0.80,0.67}
\definecolor{mediumblue}{rgb}{0.00,0.00,0.80}
\definecolor{mediumorchid}{rgb}{0.73,0.33,0.83}
\definecolor{mediumpurple}{rgb}{0.58,0.44,0.86}
\definecolor{mediumsea}{rgb}{0.24,0.70,0.44}
\definecolor{mediumslate}{rgb}{0.48,0.41,0.93}
\definecolor{mediumspring}{rgb}{0.00,0.98,0.60}
\definecolor{mediumturquoise}{rgb}{0.28,0.82,0.80}
\definecolor{mediumviolet}{rgb}{0.78,0.08,0.52}
\definecolor{midnightblue}{rgb}{0.10,0.10,0.44}
\definecolor{mintcream}{rgb}{0.96,1.00,0.98}
\definecolor{mistyrose}{rgb}{1.00,0.89,0.88}
\definecolor{moccasin}{rgb}{1.00,0.89,0.71}
\definecolor{navajowhite}{rgb}{1.00,0.87,0.68}
\definecolor{navyblue}{rgb}{0.00,0.00,0.50}
\definecolor{navy}{rgb}{0.00,0.00,0.50}
\definecolor{oldlace}{rgb}{0.99,0.96,0.90}
\definecolor{olivedrab}{rgb}{0.42,0.56,0.14}
\definecolor{orange1}{rgb}{1.00,0.65,0.00}
\definecolor{orange2}{rgb}{0.93,0.60,0.00}
\definecolor{orange3}{rgb}{0.80,0.52,0.00}
\definecolor{orange4}{rgb}{0.55,0.35,0.00}
\definecolor{orangered}{rgb}{1.00,0.27,0.00}
\definecolor{orange}{rgb}{1.00,0.65,0.00}
\definecolor{orchid1}{rgb}{1.00,0.51,0.98}
\definecolor{orchid2}{rgb}{0.93,0.48,0.91}
\definecolor{orchid3}{rgb}{0.80,0.41,0.79}
\definecolor{orchid4}{rgb}{0.55,0.28,0.54}
\definecolor{orchid}{rgb}{0.85,0.44,0.84}
\definecolor{palegoldenrod}{rgb}{0.93,0.91,0.67}
\definecolor{palegreen}{rgb}{0.60,0.98,0.60}
\definecolor{paleturquoise}{rgb}{0.69,0.93,0.93}
\definecolor{paleviolet}{rgb}{0.86,0.44,0.58}
\definecolor{papayawhip}{rgb}{1.00,0.94,0.84}
\definecolor{peachpuff}{rgb}{1.00,0.85,0.73}
\definecolor{peru}{rgb}{0.80,0.52,0.25}
\definecolor{pink1}{rgb}{1.00,0.71,0.77}
\definecolor{pink2}{rgb}{0.93,0.66,0.72}
\definecolor{pink3}{rgb}{0.80,0.57,0.62}
\definecolor{pink4}{rgb}{0.55,0.39,0.42}
\definecolor{pink}{rgb}{1.00,0.75,0.80}
\definecolor{plum1}{rgb}{1.00,0.73,1.00}
\definecolor{plum2}{rgb}{0.93,0.68,0.93}
\definecolor{plum3}{rgb}{0.80,0.59,0.80}
\definecolor{plum4}{rgb}{0.55,0.40,0.55}
\definecolor{plum}{rgb}{0.87,0.63,0.87}
\definecolor{powderblue}{rgb}{0.69,0.88,0.90}
\definecolor{purple1}{rgb}{0.61,0.19,1.00}
\definecolor{purple2}{rgb}{0.57,0.17,0.93}
\definecolor{purple3}{rgb}{0.49,0.15,0.80}
\definecolor{purple4}{rgb}{0.33,0.10,0.55}
\definecolor{purple}{rgb}{0.63,0.13,0.94}
\definecolor{red1}{rgb}{1.00,0.00,0.00}
\definecolor{red2}{rgb}{0.93,0.00,0.00}
\definecolor{red3}{rgb}{0.80,0.00,0.00}
\definecolor{red4}{rgb}{0.55,0.00,0.00}
\definecolor{red}{rgb}{1.00,0.00,0.00}
\definecolor{rosybrown}{rgb}{0.74,0.56,0.56}
\definecolor{royalblue}{rgb}{0.25,0.41,0.88}
\definecolor{saddlebrown}{rgb}{0.55,0.27,0.07}
\definecolor{salmon1}{rgb}{1.00,0.55,0.41}
\definecolor{salmon2}{rgb}{0.93,0.51,0.38}
\definecolor{salmon3}{rgb}{0.80,0.44,0.33}
\definecolor{salmon4}{rgb}{0.55,0.30,0.22}
\definecolor{salmon}{rgb}{0.98,0.50,0.45}
\definecolor{sandybrown}{rgb}{0.96,0.64,0.38}
\definecolor{seagreen}{rgb}{0.18,0.55,0.34}
\definecolor{seashell1}{rgb}{1.00,0.96,0.93}
\definecolor{seashell2}{rgb}{0.93,0.90,0.87}
\definecolor{seashell3}{rgb}{0.80,0.77,0.75}
\definecolor{seashell4}{rgb}{0.55,0.53,0.51}
\definecolor{seashell}{rgb}{1.00,0.96,0.93}
\definecolor{sienna1}{rgb}{1.00,0.51,0.28}
\definecolor{sienna2}{rgb}{0.93,0.47,0.26}
\definecolor{sienna3}{rgb}{0.80,0.41,0.22}
\definecolor{sienna4}{rgb}{0.55,0.28,0.15}
\definecolor{sienna}{rgb}{0.63,0.32,0.18}
\definecolor{skyblue}{rgb}{0.53,0.81,0.92}
\definecolor{slateblue}{rgb}{0.42,0.35,0.80}
\definecolor{slategray}{rgb}{0.44,0.50,0.56}
\definecolor{slategrey}{rgb}{0.44,0.50,0.56}
\definecolor{snow1}{rgb}{1.00,0.98,0.98}
\definecolor{snow2}{rgb}{0.93,0.91,0.91}
\definecolor{snow3}{rgb}{0.80,0.79,0.79}
\definecolor{snow4}{rgb}{0.55,0.54,0.54}
\definecolor{snow}{rgb}{1.00,0.98,0.98}
\definecolor{springgreen}{rgb}{0.00,1.00,0.50}
\definecolor{steelblue}{rgb}{0.27,0.51,0.71}
\definecolor{tan1}{rgb}{1.00,0.65,0.31}
\definecolor{tan2}{rgb}{0.93,0.60,0.29}
\definecolor{tan3}{rgb}{0.80,0.52,0.25}
\definecolor{tan4}{rgb}{0.55,0.35,0.17}
\definecolor{tan}{rgb}{0.82,0.71,0.55}
\definecolor{thistle1}{rgb}{1.00,0.88,1.00}
\definecolor{thistle2}{rgb}{0.93,0.82,0.93}
\definecolor{thistle3}{rgb}{0.80,0.71,0.80}
\definecolor{thistle4}{rgb}{0.55,0.48,0.55}
\definecolor{thistle}{rgb}{0.85,0.75,0.85}
\definecolor{tomato1}{rgb}{1.00,0.39,0.28}
\definecolor{tomato2}{rgb}{0.93,0.36,0.26}
\definecolor{tomato3}{rgb}{0.80,0.31,0.22}
\definecolor{tomato4}{rgb}{0.55,0.21,0.15}
\definecolor{tomato}{rgb}{1.00,0.39,0.28}
\definecolor{turquoise1}{rgb}{0.00,0.96,1.00}
\definecolor{turquoise2}{rgb}{0.00,0.90,0.93}
\definecolor{turquoise3}{rgb}{0.00,0.77,0.80}
\definecolor{turquoise4}{rgb}{0.00,0.53,0.55}
\definecolor{turquoise}{rgb}{0.25,0.88,0.82}
\definecolor{violetred}{rgb}{0.82,0.13,0.56}
\definecolor{violet}{rgb}{0.93,0.51,0.93}
\definecolor{wheat1}{rgb}{1.00,0.91,0.73}
\definecolor{wheat2}{rgb}{0.93,0.85,0.68}
\definecolor{wheat3}{rgb}{0.80,0.73,0.59}
\definecolor{wheat4}{rgb}{0.55,0.49,0.40}
\definecolor{wheat}{rgb}{0.96,0.87,0.70}
\definecolor{whitesmoke}{rgb}{0.96,0.96,0.96}
\definecolor{white}{rgb}{1.00,1.00,1.00}
\definecolor{yellow1}{rgb}{1.00,1.00,0.00}
\definecolor{yellow2}{rgb}{0.93,0.93,0.00}
\definecolor{yellow3}{rgb}{0.80,0.80,0.00}
\definecolor{yellow4}{rgb}{0.55,0.55,0.00}
\definecolor{yellowgreen}{rgb}{0.60,0.80,0.20}
\definecolor{yellow}{rgb}{1.00,1.00,0.00}

\usepackage[toc,page]{appendix}
\usepackage{amsmath, amssymb,paralist,amsthm}
\usepackage{float,graphicx}
\usepackage{graphics}

\def\ep{\epsilon}

\def\E{\mathbb{E}}
\def\R{\mathbf{R}}

\newtheorem{theorem}{Theorem}[section]
\newtheorem{lemma}{Lemma}[section]

\newtheorem{remark}{Remark}[section]

\begin{document}

\centerline{\Large \bf Inspection games in a mean field setting}

\bigskip
\centerline{\bf Vassili Kolokoltsov}

\centerline{\it Department of Statistics, University of Warwick}
\centerline{\it Coventry, CV4 7AL, UK, v.kolokoltsov@warwick.ac.uk}

\bigskip
\centerline{\bf   Wei Yang}

\centerline{\it Department of Mathematics and Statistics, University of Strathclyde}
\centerline{\it Glasgow, G1 1XH, UK, w.yang@strath.ac.uk}

\bigskip
\bigskip

In this paper, we present a new development of inspection games in a mean field setting. In our dynamic version of an inspection game, there is  one  inspector and a large number $N$  interacting inspectees with a finite state space.
By applying the mean field game methodology, we present a solution as an $\epsilon$-equilibrium to this type of inspection games, where $\epsilon$ goes to $0$ as $N$ tends to infinity. In order to facilitate numerical analysis of this new type inspection game, we conduct an approximation analysis, that is we approximate the optimal Lipschitz  continuous switching strategies by smooth switching strategies. We show that any approximating  smooth switching strategy is also an $\epsilon$-equilibrium solution to the inspection game with a large and finite number $N$ of inspectees with $\epsilon$ being of order  $1/N$.

\newpage

\section{Introduction}

An inspection game is a non-cooperative game whose players are often called an inspector and an inspectee. It models a situation where the inspectee, which may be an individual, an organisation, a state or a country, is obliged to follow certain regulations but has an incentive to violate them. The inspector tries to minimise the impact of such violations by means of inspections that uncover them.

A simple example of an inspection game can be described by the following 2x2 normal form game in Table \ref{The canonical two-player inspection game}, where the inspectee is the row player and the inspector is the column player, and the left (resp. right) entry of each cell corresponds to the payoffs for the individual (resp. inspector).
\begin{table}[ht]
\centering
\begin{tabular}{|l|c|c|c|c|c|}
\hline
 & Inspect & Not Inspect \\[.4em]
\hline
Violate & $-1,1$ & $2,-2$ \\[.4em]
\hline
Comply & $0,-1$ & $0,0$ \\[.4em]
\hline
\end{tabular}
\caption{The simplest two-player inspection game}
\label{The canonical two-player inspection game}
\end{table}
Typically an inspection game has a mixed equilibrium.

Inspection games were introduced by Dresher (1962) and Kuhn(1963), the underlying motivation was the cold war between the US and the Soviet Union and the desire to monitor the various arms control agreements that were signed by the two superpowers. Analytically, these settings led quite naturally to two-person game formulations with various assumptions about the strategy sets that were feasible to each of the two parties.

Inspection games have been investigated quite extensively during the last five decades. They have a wide variety of applications to name a few such as arms control by  Avenhaus et al. (1996),
auditing of accounts by Borch (1990),
tax inspection by Greenberg (1984) and Alm and McKee (2004),
environmental protection by Avenhaus (1994),
quality control in supply chains by  Reyniers and Tapiero (1995), Tapiero and Kogan (2007), Hsieh and Liu (2010),
stock keeping  by Fandel and Trockel (2008)
and communication infrastructures by Gianini et al. (2013), Chung, Hollinger and Isler (2011).
The research on inspection games contributes to the construction of an effective inspection policy for the inspector when an illegal action is executed strategically.

In the literature, attentions are mainly on two-person zero-sum games, with some drifts to two-person non-zero-sums models. As far as we know, there are very few models with multiple inspectees. In the arms control inspection context, Kilgour and Averhaus (1994) considered a model with two or more inspectees, where the inspectees are independent and the inspector's inspection is a binary variable, namely inspector decides whether to inspect or not.  Later Avenhaus and Kilgour (2004) studied another model with two inspectees, where the inspector has a fixed level and continuously divisible inspection resources. They answered the question on how the inspector distributes efficiently its limited inspection resources over several independent inspectees.

In a recent work on distributed information systems, Gianini et al. (2013) study a simultaneous one-shot  inspection game with uncoordinated $m$ inspectors and $n$ non-interacting inspectees. In their model, an inspector has limited resources but the probability of detection is not a function of the resourcee.
They show that due to the lack of coupling among inspectees, adding or removing inspectees does not change the best mixed strategy of one inspectee.

Kolokoltsov, Passi and Yang (2013) develop a dynamic inspection game model with one inspector and a large number of  interacting inspectees in an evolutionary setting. In that model, each inspectee is under evolutionary pressure and periodically updates their strategies after binary interactions.
Specifically, at the beginning of each period, an exogenous and fixed fraction of the population can update their behaviour upon meeting another randomly chosen individual in the population. If two inspectees meet and have the same strategy, then that strategy is retained by the updating individual. If however, the two individuals have different strategies, the updating individual may revise his behaviour on the basis of the payoffs enjoyed by the two in the previous period. Solutions in terms of mixed strategies to this type of games are presented therein.

In this paper, we present a new development of inspection games in a mean field setting. In this new model of inspection games,
there is one inspector and a large number of interacting inspectees. Different from the settings by Kilgour and Avenhaus (1994),  Avenhaus and Kilgour (2004) and Gianini et al. (2013) where the multiple inspectees are independent, our new model considers interacting inspectees in the sense that one inspectee's payoff does depend on one another's strategy.  Also different from the evolutionary setting by Kolokoltsov, Passi and Yang (2013) where the binary interactions are considered and inspectees have no payoff functions, in our new model each inspectee aims to maximise her own total expected payoff which depends on the aggregate behaviour of the population.
We aim to approximate the Nash solution to this new type of dynamic inspection games by using the recently developed theory of mean field games.

The mean field game theory is a new branch of game theory and it has become a powerful tool to study complex games with a large number of players. The initial work done by Lasry and Lions  (2006a, 2006b, 2007) and Huang, Malham\'e and Caines (2006, 2007) consider continuous-time continuous-state games. In this paper in the context of inspection games, we study a mean field game in a finite state space setting. The finite state space setting is also considered within the context of socio-economic sciences by Gomes, Mohr and Souza (2010, 2013), Gomes, Velho, Wolfram (2014a, 2014b), wherein the authors study one hyperbolic equation and focus on the analysis of shock-formation phenomenon of the system in the setting of  two-state mean field games. Though those above mentioned papers consider a similar framework as the one in this paper, different approaches are applied.
In this paper, we apply a modified version of the standard mean field games method, namely we study a system of coupled two differential equations with one being forward and the other one backward. Our model contributes to the study of mean field games on a finite state space with a major player.

The paper is organised as follows.
Section 2 describes in detail the model of an inspection game with one inspector and a large number $N$ inspectees.
In Section 3, we set up a mean field inspection game with a continuum of inspectees and derive the system of coupled equations  \eqref{limit HJB}-\eqref{limit distribution}.
Main results of this paper are shown in Section 4.
First, in Theorem \ref{Thm1}, we show that for a short time game, there exists a unique solution to the mean field inspection game, namely, the single inspector has a unique best response to the continuum of inspectees and any representative inspectee has a unique best response to the inspector and the aggregate behaviour of  the continuum; whereas, for a long time game, we show the existence of a solution to the mean field inspection game. Then in Theorem \ref{convergence of propagators}, we show that the probability distributions of finite $N$-inspectees on their state space converges to the one of a continuum limit as $N\to\infty$.  Finally, in Theorem \ref{e convergence} we conclude that an optimal Lipschitz continuous switching strategy, which is derived from the mean field inspection game, is an $\epsilon$-equilibrium to an inspection game with a finite number $N$ of inspectees, with $\epsilon=\epsilon(N)\to 0$ as $N\to\infty$.
In Section 5 we conduct an approximation analysis in order to adapt our theoretical results for numerical analysis and to discuss the  rate of convergence.  We approximate optimal Lipschitz continuous switching strategies $q^*$ by a sequence of smooth switching strategies $q^*_\eta$, $\eta>0$. We show that any approximating  smooth switching strategies $q^*_\eta$ is also an $\epsilon$-equilibrium solution to the finite $N$ inspection game, with $\epsilon=\epsilon(N,\eta)\to 0$ as $N\to \infty$ and $\eta\to 0$. Further, using smooth switching strategies $q^*_\eta$ as an $\epsilon$-equilibrium solution, we show that $\epsilon=\epsilon(N,\eta)$ is of order  $1/N$.

\section{An inspection game with $N$ interacting inspectees}
\numberwithin{equation}{section}

We consider a dynamic inspection game in a continuous time setting with a finite time horizon $T>0$. In this game, there is one inspector and $N$ (a fixed integer) inspectees (refereed to as the population). Roughly speaking, every inspectee chooses her crime levels they would commit to maximise her payoff function. The controlled dynamics of crime levels for each inspectee is modelled as a controlled continuous-time Markov Chain.  The inspector decides the amount of investment for inspection so as to maximise her payoff, based on the observation of the crime distribution. To study this game, normally one looks for a profile of best responses for all inspectees and the inspector.

Formally, first we  discuss the $N$ inspectees.  Let $\mathbb{L}_d=\{l_1,\dots,l_d\}$, $d\in\mathbb{N}$, be the state space of any inspectee. States $l_i\in \mathbb{L}_d$, $i=1,\dots, d$, are interpreted as crime levels; in other words, $l_i$ can be understood as illegal profits one can gain by committing crimes.

Denote by $\Sigma_d$ the set of probability distributions on $\mathbb{L}_d$, i.e.
\begin{equation}
\label{Limit space space}
\Sigma_d=\{x=(x_1,\dots, x_d)\in[0,1]^d: \sum_{j=1}^d x_j=1\}
\end{equation}
and by $C([0,T], \Sigma_d)$ the set of continuous curves $\{x(t)\in \Sigma_d, t\in[0,T]\}$, equipped with the norm
\begin{equation}
\label{norm}
\|x(\cdot)\|_{\infty}=\sup_{t\in[0,T]}\|x(t)\|
\end{equation}
where $\|\cdot\|$ denotes the  Euclidean norm in $\R^d$.

The dynamics of  every inspectee is modelled by a continuous-time Markov chain on $\mathbb{L}_d$. Every inspectee chooses a switching strategy between the crime levels to maximise her own objective function. Specifically, the dynamics of  the crime levels of the inspectee $a\in\{1,\dots, N\}$ is modelled by a continuous-time Markov chain $M^{(a)}=\{M^{(a)}(t), t\in[0,T]\}$ on the state space $\mathbb{L}_d$. For a given curve $\{x(t), t\in[0,T]\}\in C([0,T], \Sigma_d)$, the stochastic dynamic $M^{(a)}$ is specified by  the switching matrix 
\begin{align}
\label{switching matrix}
\mathbb{Q}^{a}(t,x(t))= \left( \begin{array}{ccc}
 & q^{a}(t,l_1,x(t))&  \\
 & \vdots &  \\
& q^{a}(t,l_i,x(t))&  \\
 & \vdots &  \\
& q^{a}(t,l_d,x(t))& 
 \end{array} \right)=
\left( \begin{array}{ccc}
q^{a}_{11}(t,x(t)) & \dots & q^{a}_{1d}(t,x(t)) \\
 \vdots &  &  \vdots \\
 q^{a}_{i1}(t,x(t)) & \dots & q^{a}_{id}(t,x(t)) \\
 \vdots &  &  \vdots \\
q^{a}_{d1}(t,x(t)) & \dots & q^{a}_{dd}(t,x(t))
 \end{array} \right)
 \end{align}
which is chosen by the $a$th inspectee. At any time $t\in[0,T]$, $x\in\Sigma_d$ and $j\neq i$, the entry $q^a_{ij}(t,x)$ is in $[0,Q]$ and presents the infinitesimal transition rate from state $l_i$ to state $l_j$, bounded by a  constant $Q>0$. Moreover for any $i$, $\sum_{j=1}^d q^a_{ij}(t,x)=0$, namely $q^a_{ii}(t,x)$ is chosen in such a way that  $q^a_{ii}(t,x)=-\sum_{j\neq i} q^a_{ij}(t,x)$.



We are interested in symmetric inspectees, meaning that any inspectees who are at the same time at the same crime level will choose the same switching rate. In other words, the choice of a switching strategy does not depend on the identity  of each inspetee. Thus we can omit the identity index and for any $(t,l_i,x)\in [0,T]\times \mathbb{L}_d\times \Sigma_d$ denote
$$q(t,l_i, x): =q^{(a)}(t,l_i, x)\quad \text{and } \quad q_{ij}(t,x)=q^a_{ij}(t,x).$$

Next we discuss the dynamics of the population which is  denoted by 
$$X^N=\left\{\big(X^N_1(t), \dots,X^N_d(t)\big), t\in[0,T]\right\}.$$
The dynamics $X^N$ is a Markov process and describes the evolution of crime distributions, that is for any $i=1,\dots,d$
$$X^N_i(t)=\frac{\sharp\{a\in\{1, \dots, N\}: M^{(a)}(t)=l_i\}}{N}$$
specifies the proportion of inspectees at the crime level $l_i$ at time $t\in[0,T]$.
The superscript $N$  in $X^N$ is used to distinguish between the dynamics in the finite population of size $N$ and the one in the limit (to be introduced in Section 3). The state space of the population is denoted by
$$\mathbb{S}_d^{{N}}= \left\{y=\left (\frac{n_1}{N}, \dots, \frac{n_d}{N}\right): \sum_{j=1}^dn_j=N\right\},$$
which is a subset of the closed simplex $\Sigma_d$ defined in \eqref{Limit space space}.

%
%

The Markov process $X^N$ on the state space $\mathbb{S}_d^{{N}}$  is generated by the time-inhomogenous operator 
$L_t^{{N}}:C(\mathbb{S}_d^{{N}})\to C(\mathbb{S}_d^{{N}})$
defined by
\begin{equation}
\label{L(1)}
L_{t}^{{N}}f(y)=\sum_{\substack{i,j=1 \\ i\neq j}}^d (y\cdot e_i) q_{ij}(t,y) N \left [   f\left(y-\frac{e_i}{N}+\frac{e_j}{N}\right)-f(y)\right],
\end{equation}
where $e_i$, $i\in\{1,\dots,d\}$, denotes the standard basis in $\R^d$, namely the $i$th entry is $1$ and all the other entries are $0$; $(y\cdot e_i)=n_i/N$ for any $y=(n_1/N, \dots, n_d/N)$. 


An intuitive probabilistic interpretation of the stochastic process $X^N=\{(X^N_1(t), \dots,X^N_d(t)), t\in[0,T]\}$ is as follows.  At the initial stage of this game $t=0$, the $N$ inspectees are distributed  arbitrarily among the $d$ states $\{l_1,\dots,l_d\}$. 
The initial state of the population is described by the vector $X^N(0)=(n_1(0)/N,\dots, n_d(0)/N)$. Here  $n_j(0)$, $j\in\{1,\dots,d\}$, specifies the number of inspectees at the crime level $l_j$ at $t=0$.   As the dynamic of each inspectee is modelled as a Markov chain with a switching matrix in the form of \eqref{switching matrix}, every inspectee has a random waiting time at their current crime level before she switches to another crime level. 
Denote by $\tau_1$ the shortest waiting time among $N$ inspectees. Then at the time $\tau_1$, it is the first time when an inspectee changes her crime level, say from $l_i$ to $l_j$. Consequently, from the initial state $X^N(0)$, the Markov process $X^N$ obtains the new state:
$$X^N(\tau_1)= X^N(0)-\frac{e_i}{N}+\frac{e_j}{N}= \left(\frac{n_1(0)}{N},\dots, \frac{n_i(0)-1}{N}, \dots, \frac{n_j(0)+1}{N}, \dots, \frac{n_d(0)}{N}\right).$$
Then the process $X^N$ evolves in the manner as described above from the new state $X^N(\tau_1)$.

It is important to note that the curve $\{x(t), t\in[0,T]\}$ in \eqref{switching matrix} is the realisation of the crime distribution evolution $\{(X^N_1(t), \dots,X^N_d(t)), t\in[0,T]\}$ in the population of $N$ inspectees  or $\{(X_1(t), \dots,X_d(t)), t\in[0,T]\}$ in the continuum of inspectees (see in the following Eq. \eqref{limit distribution} ).
This is the exact place to see how inspectees interact with each other: the state dynamics of $a$th insepctee $M^{(a)}$ is influenced by all other inspectees strategies through the crime distribution; in other words, any inspectee's state dynamics is influenced by the aggregated behavior of all inspectees. This is what we mean by the {\it mean field interaction setting}.


In this model, we do not require that every inspectee has perfect information on the crime levels of all other inspectees, but we assume that everyone has the access to the exact information about the aggregate behaviour of the population, namely the crime distribution. 

Next, we introduce the single inspector.
The inspector has limited available resources for inspection, denoted by $F>0$. At any time $t\in[0,T]$, she needs to decide an  amount of resource $\alpha (t)\in[0, F]$ to be invested in inspection in order to maximise her total expected payoff function. It is assumed that the inspection resource is uniformly distributed among the population and the inspector will charge a fine $\sigma l_j$, $\sigma>0$, if she uncovers a crime at the level $l_j\in \mathbb{L}_d$.

Due to the limited resources, a complete surveillance of all inspectees' actions is practically not possible. Therefore, inspection takes place in form of a randomisation. We introduce a  detection function, which relates the detection probability to inspection resources and is given by $P: [0,F]\to [0,1]$. It is assumed that the detection function $P$ is nondecreasing and concave in inspection resources $\alpha$, i.e.
\begin{equation}
\label{detection function}
P'(\alpha)>0\quad \text{and  } \,P''(\alpha)<0,\quad \text{for } \alpha \in[0,F].
\end{equation}

One way to understand  the value $P(\alpha)$ for an $\alpha \in[0,F]$ is that, every inspectee is inspected and an illegal behaviour can be uncovered with the probability $P(\alpha)$. Another interpretation of $P(\alpha)$ is that, with the inspection resources $\alpha$ invested, a proportion $P(\alpha)$ of the population is inspected with perfect inspection, in other words, every inspectee has a probability $P(\alpha)$ to be inspected and an illegal behaviour will be detected with probability $1$.

Now, we are ready to  introduce the payoff functions for inspectees and the single inspector.
At each time $s\in[0,T]$ with a crime distribution $X^N(s)$ of the population, if the inspector invests $\alpha(s)$ for inspection, the $a$th inspectee with her crime level $M_N^{(a)}(s)$ faces a probability $P(\alpha (s))$ with which her illegal behaviour will be detected and the inspectee has to pay a fine $\sigma M_N^{(a)}(s)$; on the other hand,  she escapes with a probability $1-P(\alpha (s))$ and gains an illegal profit $M_N^{(a)}(s)$. Moreover, the inspectees pays a cost for the change of strategies, which is quadratic in the transition rates
i.e., $\sum_{l_j\neq M_N^{(a)}(s)} q^2_{M_N^{(a)}(s) \,l_j}(s, X^N(s))$.
Therefore, at each time $s$ with a crime distribution $X^N(s)$ of the population, the $a$th inspectee with a crime level $M_N^{(a)}(s)$ has a running payoff
$$   (1-P(\alpha(s))) M_N^{(a)}(s) - P(\alpha(s))\sigma M_N^{(a)}(s)-\sum_{l_j\neq M_N^{(a)}(s)} q^2_{M_N^{(a)}(s)  \,l_j}(s, X^N(s) )$$
and
a terminal payoff $J_T(M_N^{(a)}(T), X^N(T))$, where $J_T:\mathbb{L}_d\times \Sigma_d\to \R$. Therefore,
for a given $X^N\in C([0,T], \Sigma_d)$,  the $a$th inspectee aims to mamixise her  payoff  function
\begin{align}
\label{insepctee payoff}
J^{(a)}(t,l_i, q_i; X^N)=\E_{l_i}\Big[\int _t^T \Big [ &(1-P(\alpha(s))) M_N^{(a)}(s) - P(\alpha(s))\sigma M_N^{(a)}(s)\notag\\
&-\sum_{l_j\neq M^{(a)}(s)} q^2_{M_N^{(a)}(s) \, l_j}(s, X^N(s))\Big]ds+J_T(M_N^{(a)}(T), X^N(T))\Big]
\end{align}
over switching strategies $q_i=\{(q_{i1}(t), \dots, q_{id}(t)), t\in[0,T]\}$ with $q_{ij}(t)\in[0,Q]$ for any $t\in[0,T]$ and $j=1,\dots, d$, $\sum_{j=1}^d q_{ij}(t)=0$, and $q_{ii}(t)$ is such that  $q_{ii}(t)=-\sum_{j\neq i} q_{ij}(t)$.

In the meanwhile, the inspector can observe the crime distribution $X^N(s)$ of the population. At any time $s\in[0,T]$,  the inspector pays an amount of investment resources for inspection $\alpha(s)$ and gets a payoff from any individual inspectee at the crime level $l_i$:
\begin{equation}\label{single contribution}
\Phi^N_i(s):=\eta^N\Big(P(\alpha(s)) \sigma l_i - (1-P(\alpha(s))) l_i\Big),
\end{equation}
where  $\eta^N>0$ is given and prescribes the weight that the inspector assigns to any single inspectee in the finite $N$ population. The inspector aims to maximise her expected payoff at each time instance $s\in[0,T]$, that is she wants to maximise
\begin{equation}
\label{Inspector-payoff}
\E \left (-\alpha(s)+\sum_{i=1}^d N  X_i^N(s)\Phi^N_i(s)\right)
\end{equation}
over $\alpha(s)\in[0,F]$. Plug \eqref{single contribution} into \eqref{Inspector-payoff} and define {\it $L:=N\eta^N$}. Then the inspector aims to maximise her payoff function 
\begin{equation}
\label{Inspector-payoff-1}
U_N(\alpha (s), X^N(s))=\E \left (-\alpha(s)+L \sum_{i=1}^d   X_i^N(s)\Big(P(\alpha(s)) \sigma l_i - (1-P(\alpha(s))) l_i\Big)\right).
\end{equation}
over her inspection investment $\alpha(s)$.

By differentiating the function $U_N$ with respect to the first variable and together with Condition \eqref{detection function}, the function $U_N$ in Eq. \eqref{Inspector-payoff-1} has a unique maximiser $ \alpha_N^*: \Sigma_d\to [0, F]$:
\begin{align}
\label{Inspector best response}
 \alpha_N^* (s)
 = \alpha_N^* (X^N(s)):&=\arg \max_{\alpha \in[0,F]}U_N(\alpha,X^N(s))\notag\\
&=\min \left\{(P')^{-1}\left(\frac{1}{L \big(1+\sigma)\E \left[\sum_{i=1}^d  l_iX_i^N(s)\right]}\right), F\right\}
\end{align}
where $(P')^{-1}$ denotes the inverse function of $P'$. 

Note that the index $N$  is used  in the notations of objects in the setting of finite $N$ inspectees to distinguish them from their counterparties in the mean field inspection game, see Section 3.

\begin{remark}
The parameter $\eta^N$ is small, compared to the large number $N$ of inspectees. The parameter $L=N\eta^N$ in \eqref{Inspector-payoff-1} can be interpreted as the approximate total fine from the whole population of the inspectees. In Section
 3, we will consider a limiting model by sending the number of inspectees to infinity, i.e. $N\to \infty$. In the limiting model, the value of the parameter $L$ is kept the same as in the finite population problem and $\eta^N \to 0$ as $N\to \infty$. In other words, in the inspector's viewpoint, when the number of inspectees becomes very large, her total expected payoff is always bounded, and any individual inspectee's contribution $\Phi^N$ in \eqref{single contribution} to the inspector's payoff  becomes negligible. 
\end{remark}

We are interested in finding Nash equilibria of this type of inspection games. This means to find the best-response investment strategy  $\alpha_N^*(t)$ at any  $t\in[0,T]$ and a family  $\{\mathbb Q^{*(1)}, \dots, \mathbb Q^{*(N)}\}$, where $\mathbb Q ^{*(a)}, a\in\{1,\dots, N\}$, denotes the best switching strategy of the $a$th inspectee as the best response.  When $N$ is very large, the complexity of this problem gets immense.

In this paper, we will apply the mean field games methodology to solve this type of games and provide an $\epsilon$-equiliblium. Roughly speaking, we will take the number $N$ of inspectees to infinity and set up the (limiting) mean field model with a continuum of inspectees. We call this game with a single inspector and a continuum of inspectees a {\it mean field inspection game}.  Then we prove that any solution to the mean field inspection game presents an $\epsilon$-equilibrium to the original model with $N$ inspectees.

\section{The mean field inspection game}

In this section, we study a mean field inspection game with a continuum of inspectees and one inspector.
First, let's discuss the dynamics of the continuum limit as $N\to\infty$. The state space of the continuum is naturally specified by $\Sigma _d$ in \eqref{Limit space space}. Observe that, for $f\in C^1(\Sigma _d)$,
$$ \lim_{\substack{N\to \infty\\y\to x}} N  \left [   f(y-\frac{e_i}{N}+\frac{e_j}{N})-f(y)\right]=\frac{\partial f}{\partial x_j}(x)-\frac{\partial f}{\partial x_i}(x)$$
so that the limiting generator  $A_{t}:C^1(\Sigma _d)\to C(\Sigma _d)$ of $L_{t}^{{N}}$ in \eqref{L(1)} as $N\to\infty$ is of the form
 \begin{align}
\label{limit operator}
A_{t} f(x) : &= \lim_{\substack{N\to \infty\\y\to x}}   L_{t}^{{N}}f(y)\notag\\
&=  \sum_{\substack{i,j=1 \\ i\neq j}}^dx q_{ij}(t,x) \left( \frac{\partial f}{\partial x_j}(x)-\frac{\partial f}{\partial x_i}(x) \right) \notag\\
&= \sum _{\substack{j=1 \\ i\neq j}}^d\left(x_iq_{ij}(t,x) -x_j q_{ji}(t,x)\right)  \frac{\partial f}{\partial x_j}(x).
\end{align}
The limiting operator $A_t$ in \eqref{limit operator} generates a controlled crime distribution evolution of  the continuum of inspectees on the state space $\Sigma_d$, which is denoted by $X=\{X(t)=(X_1(t), \dots, X_d(t)): t\in[0,T]\}$, where $X_i(t)$ describes the fraction of inspectees at the crime level $l_i$.  Then the limiting distribution evolution $X$ is governed by  the kinetic equation
\begin{equation}
\label{Kinetic}
\frac{dX_i(t)}{dt}= \sum _{j=1}^d X_j(t) q_{ji}(t, X(t)), \quad i=1 ,\dots, d.
\end{equation}
It is worth noting that in contrast to $X^N$ generated by $L_t^N$ in  \eqref{L(1)},  the limiting evolution $X$ is deterministic, since it is modelled as the solution to the ordinary differential equation \eqref{Kinetic}.


Now in the mean field inspection game, at any time $t\in[0,T]$ with a crime distribution $X(t)=(X_1(t),\dots,X_d(t))$, the inspector aims to maximise her payoff  $U(\alpha(t), X(t))$ over her inspection investment $\alpha(t)$, where
\begin{align}
\label{inspector payoff}
U:[0,F]\times \Sigma_d\to \R,\qquad U(\alpha,x):= -\alpha+L \big((1+\sigma) P(\alpha) -1\big)\sum_{i=1}^d  l_ix_i.
\end{align}
By the condition \eqref{detection function}, for any $X(t)\in\Sigma_d$, the function $U$ in \eqref{inspector payoff} has a unique maximiser $ \alpha^*(t)$ given by:
\begin{align}
\label{limiting-Inspector-payoff-1}
 \alpha^*(t)=\alpha^* (X(t)): &=\arg \max_{\alpha \in[0,F]}U(\alpha,X(t))\notag\\
&=\min \left\{(P')^{-1}\left(\frac{1}{L \big(1+\sigma)\sum_{i=1}^d  l_iX_i(t)}\right), F\right\}
\end{align}
where $(P')^{-1}$ denotes the inverse function of $P'$. 

Let $\{M(t), t\in[0,T]\}$ denote the  controlled state dynamics of a representative inspectee of the continuum of inspectees. Starting at any time $t\in[0,T]$ and a state $l_i\in\mathbb{L}_d$, a representative inspectee of the continuum aims to maximise
the payoff
\begin{align}
\label{insepctee payoff-2}
 \E_{l_i}\Big[\int _t^T  \Big [ (1-&P( \alpha^*(X(s))) M(s) - P( \alpha^*(X(s)))\sigma M(s)\notag\\
&- \sum_{l_j\neq M(s)} q^2_{M(s)l_j}(s,X(s))\Big]ds + J_T(M(T), X(T))\Big]
\end{align}
 where the function $\alpha^*$ is defined in \eqref{limiting-Inspector-payoff-1} and the terminal cost function $J_T: \mathbb{L}_d\times\Sigma_d\to \R$.
 Denote the running cost function in the payoff \eqref{insepctee payoff-2} by $h: \mathbb{L}_d\times \Sigma_d\times\R^d\to \R$
 \begin{align}\label{running cost}
 h(l_i,x,q_i)&=  (1-P(\alpha^*(x)))l_i -P( \alpha^*(x))\sigma l_ i -\sum_{j\neq i} q^2_{ij} \\
 &= l_i -l_i(1+\sigma) P( \alpha^*(x))-\sum_{j\neq i} q^2_{ij} \notag
 \end{align}
where $q_i=(q_{i1},\dots, q_{id})$.
The value function  for a representative inspectee $V:[0,T]\times\mathbb{L}_d\times C([0,T],\Sigma_d)\to \R$ is defined as
\begin{align}
\label{Value function}
V(t,l_i;X)=\sup _{q(\cdot,\cdot,\cdot)} \E_{t,l_i}\Big[\int _t^T &h\big(M(s), X(s), q(s,M(s), X(s))\big) ds + J_T(M(T), X(T))\Big]
\end{align}
over measurable functions $q:[t,T]\times \mathbb{L}_d\times \Sigma_d\to \R^d$. For any $t\in[0,T]$ and $X\in C([0,T],\Sigma_d)$, denote the norm of a value function $V$ on $\mathbb{L}_d$ by
 \begin{equation}
 \label{norm V}
 \|V(t,\cdot;X)\|:=\sup_{l_i\in\mathbb{L}_d} |V(t, l_i;X)|.
 \end{equation}

By the dynamic programming principle, for any given distribution evolution $X\in C([0,T],\Sigma_d) $,  the value function $V$ in \eqref{Value function} satisfies Hamilton-Jacobi-Bellman (HJB) equation
\begin{align}
\label{HJB-g}
\frac{d V}{d t}(t,l_i)+H(l_i, V(t,\cdot), X(t))=0
\end{align}
with a terminal function $V(T,\cdot)=J_T(\cdot, X(T))$,
where  the function $H: \mathbb{L}_d\times \R^d\times \Sigma_d\to\R$
\begin{align}
\label{H-def}
H(l_i,\phi,x)&=\sup_{{q_i\in\R^d}} \left[  l_i -l_i(1+\sigma) P(\alpha^*(x))-\sum_{j\neq i} q^2_{ij} +    \sum _{j\neq i} (\phi_j-\phi_i) q_{ij} \right]
\end{align}
with the function $ \alpha^*$ defined in \eqref{limiting-Inspector-payoff-1}, where $\phi=(\phi_1,\dots, \phi_d)$ and $q_i=(q_{i1}, \dots, q_{id})\in \R^d$.
The first order condition shows that the unique maximiser in \eqref{H-def} is 
\begin{equation}
\label{optimal q 1}
q^*(t,l_i,x;\phi)=(q_{i1}^*(t,x;\phi),\dots, q_{id}^*(t,x;\phi))
\end{equation}
where for $ j\neq i$
\begin{equation*}
 q^*_{ij}(t,x;\phi)=\left\{
\begin{array}{c l}
    0 & \text{if } \phi_j-\phi_i<0\\
   \frac{1}{2}(\phi_j-\phi_i) &  \text{if } 0\leq \phi_j-\phi_i\leq 2Q \\
    Q &\text{if } \phi_j-\phi_i>2Q
\end{array}
\right.
\end{equation*}
and by definition $q^*_{ii}(t,x;\phi)=- \sum _{j\neq i} q^*_{ij}(t,x;\phi)$. 

Notice that the individual optimal switching function $ q^*$ in \eqref{optimal q 1} does not explicitly depend on the aggregate behaviour of the population $X$. However, in the model, the variable $\phi$ will be the value of the value function in \eqref{Value function} at each time instance, which depends on $X$.

To summarise, the function $\alpha^*$ in \eqref{limiting-Inspector-payoff-1} gives the best response  for any inspector to the prevailing crime distribution $X(t)\in\Sigma_d$ at any time $t\in[0,T]$.
The running cost function $h$  in \eqref{running cost} for a representative inspectee depends on the  best response function of the inspector $\alpha^*$. The value function $V$ in \eqref{Value function} for the representative inspectee is the solution of  the HJB equation \eqref{HJB-g}. The resulting best response of a representative inspectee is the optimal switching function $q^*$ in \eqref{optimal q 1}.
Rationally, the representative inspectee  applies the resulting optimal control policy $q^*$ to control her state dynamics of crime levels. The resulting  crime distribution evolution, as the aggregate behaviour of the continuum of inspectees,  is described by the solution  $X=\{X(t), [0,T]\}$ to \eqref{Kinetic} with the optimal switching $q^*$ in \eqref{optimal q 1}. The solution $X$ should be consistent with the one observed by the inspector as the prevailing crime distribution. Hence we get the following system of coupled equations
\begin{align}
\label{limit HJB}
\left\{
\begin{array}{c l}&\displaystyle
\frac{d V(t,l_i)}{d t}+H(l_i,V(t,\cdot),X(t)) =0,  \quad\quad \quad\quad i=1 ,\dots, d \\[.6em]
&V(T,\cdot)=J_T(\cdot, X(T))
\end{array}
\right.
\end{align}

\begin{align}
 \label{limit distribution}
\left\{
\begin{array}{c l}
&\displaystyle\frac{dX_i(t)}{dt}=\sum _{j=1}^d X_j(t) q^*_{ji}(t,X(t);  V(t,\cdot) ), \quad i=1 ,\dots, d \\[.6em]
& X(0)=x(0)
\end{array}
\right.
\end{align}
where the Hamiltonian $H$ is defined in \eqref{H-def} and the optimal transition function $q^*$ is defined in \eqref{optimal q 1}.
The main feature of this coupled system of equations is that Eq. \eqref{limit HJB} is a backward ordinary differential equation, yet  \eqref{limit distribution}  is a forward ordinary differential equation. This model
can be viewed as a modified version of standard mean field games equations in discrete-state space setting with a deterministic major player.

\section{Main results}

First, we discuss the existence and uniqueness of a solution to the coupled system \eqref{limit HJB}-\eqref{limit distribution}. In Theorem \ref{Thm1}, we show that for a short time game, there exists a unique solution to the mean field inspection game, namely, the single inspector has a unique best response to the continuum of inspectees and any representative inspectee has a unique best response to the inspector and the aggregate behaviour of  the continuum; whereas, for a long time game, we show the existence of a solution to the mean field inspection game. Then in Theorem \ref{convergence of propagators}, we show that the probability distributions of finite $N$-inspectees on their state space converges to the one of a limiting system as $N\to\infty$.  Finally, in Theorem \ref{e convergence}, we conclude that a solution to the mean field inspection game is an $\epsilon$-equilibrium to an inspection game with a finite number $N$ of inspectees.

\begin{theorem}
\label{Thm1}

(i) For a small $T$, the coupled system of equations \eqref{limit HJB}-\eqref{limit distribution} has a unique solution
  $(X, V)$;

 (ii) for an arbitrary finite $T$, there exists a solution to the coupled system of equations  \eqref{limit HJB}-\eqref{limit distribution}.
\end{theorem}
\proof
 The proof to Theorem \ref{Thm1}  consists of three steps and can be found in Appendix \ref{Proof to Theorem Thm1}
 \qed

Next, we prove that the pair of resulting best responses $\alpha^*(X(t))$  and $q^*(t,l_i, X(t); V(t,\cdot ;X))$ presents an $\epsilon$-equilibrium of an inspection  game with one inspector and $N$ inspectees, where $(X,V)$ is a solution to equations  \eqref{limit HJB}-\eqref{limit distribution} and $\alpha^*$ and $q^*$ are defined respectively in \eqref{limiting-Inspector-payoff-1}  and \eqref {optimal q 1}.

To this end, we will tag one insepctee and impose that she  applies a switching strategy $\tilde q (t,l_i,x)$ which is Lipschitz continuous in the variable $x$ and different from $q^*(t,l_i,x)$. Let $\tilde M^{tag, \tilde q}_N(t)$ (resp. $\tilde M^{tag, q^*}_N(t)$) denote the  state dynamics of the tagged inspectee with switching strategy $\tilde q$ (resp. with $q^*$) in the finite $N$ inspectees setting with a given initial data $m_N(0)\in\mathbb{L}_d$.
Let $\tilde M^{tag, \tilde q}(t)$ (resp. $\tilde M^{tag, q^*}(t)$) denote the  state dynamics of the tagged inspectee with switching strategy $\tilde q$ (resp. with $q^*$) in the continuum limit  with  a given  initial data $m(0)\in\mathbb{L}_d$. Meanwhile all other inspectees apply the same strategy $q^*$.

The controlled Markov process $\{\big(X_{[q^*, \tilde q]}^N(t), \tilde M^{tag, \tilde q}_N(t)\big), t\in[0,T]\}$ of $N$ interacting inspectees  is generated by the operators $\widehat L_t ^N[q^*, \tilde q]$ acting on $[C(\mathbb{S}_d^{{N}}\times \mathbb{L}_d)]^d$ (the set of continuous and bounded vector-valued functions $f$ on $\mathbb{S}_d^{{N}}\times\mathbb{L}_d$):
\begin{align}
\label{tagged L}
\widehat L_t ^N[q^*, \tilde q] f(y,l_k)=&\sum_{\substack{i,j=1 \\ i\neq k}}^d (y\cdot e_i) q^*_{ij}(t,y) N \left [   f(y-\tfrac{e_i}{N}+\tfrac{e_j}{N}, l_k)-f(y,l_k)\right]\notag\\
&+ \sum_{j=1}^d (y\cdot e_k-\tfrac{1}{N}) q^*_{kj}(t,y) N \left [   f(y-\tfrac{e_{k}}{N}+\tfrac{e_j}{N}, l_k)-f(y,l_k)\right]\notag\\
&+\sum_{j=1}^d \tfrac{1}{N} \tilde q_{kj}(t,y) N \left [   f(y-\tfrac{e_k}{N}+\tfrac{e_j}{N}, l_j)-f(y,l_k)\right]
\end{align}
where $y=(n_1/N,\dots,n_d/N)$.

One can under the operators in \eqref{tagged L} in the following way. Consider that the tagged inspectee is at the crime level $l_k\in\mathbb{L}_d$ and there are $n_k$ inspectees who are at the crime level $l_k$.
The first term in \eqref{tagged L} prescribes the interactions between inspectees who are any crime levels $l_i\in\mathbb{L}_d\setminus l_k$, and apply the strategy $q^*$; the second term in \eqref{tagged L} prescribes the interactions between inspectees who are at the crime level $l_k$, except the tagged inspectee, and apply the strategy $q^*$; the third term in \eqref{tagged L} prescribes the behaviour of the tagged inspectee who are at $l_k$ and applies the strategy $ \tilde q$.

\begin{remark}
In the definition of the operator $\widehat L_t ^N[q^*, \tilde q]$ in \eqref{tagged L}, the position of the tagged inspectee is counted twice as his own position (crime level) $\tilde M^{tag, \tilde q}_N(t)$ and as his contribution towards the crime empirical measure $X_{tag}^N(t)$. Alternatively, instead of $X_{[q^*, \tilde q]}^N(t)$ one can take the empirical measure of other inspectees only, i.e. $X_{[q^*]}^{N-1}(t)$. Such notations are used e.g. by Gomes in \cite{GVW2014-2}. However, the distinction between these two notations disappear in the limit $N\to\infty$.
\end{remark}

The controlled process $\{\big(X_{[q^*, \tilde q]}(t), \tilde M^{tag, \tilde q}(t)\big), t\in[0,T]\}$ of the continuum limit as  $N\to \infty$ is generated by the limiting operators  $\widehat A_t[q^*, \tilde q]: [C^1(\Sigma_d\times \mathbb{L}_d)]^d\to [C(\Sigma_d\times\mathbb{L}_d)]^d$ of $\widehat L_t ^N$ in \eqref{tagged L}:
\begin{align}
\label{tagged A}
\widehat A_t[q^*, \tilde q]f(x,l_k):=& \lim_{\substack{N\to \infty\\y\to x}}   \widehat L_{t}^{{N}}f(y,l_k)\notag\\
=& \sum _{\substack{j=1 \\ i\neq j}}^d\left(x_iq^*_{ij}(t,x) -x_j q^*_{ji}(t,x)\right)  \frac{\partial f}{\partial x_j }(x,l_k)+ \sum _{\substack{j=1 \\ i\neq j}}^d \tilde q_{kj(t,x)}(f(x,l_j)-f(x,l_k))
\end{align}
where the space $[C^1(\Sigma_d\times \mathbb{L}_d)]^d$ is the set of continuous and bounded vector-valued functions $f$ on $\Sigma_d\times\mathbb{L}_d$ which are differentiable in the first variable.

\begin{remark}\label{Markov processes}
It is worth noting that the process $\{\big(X_{[q^*, \tilde q]}^N(t), \tilde M^{tag, \tilde q}_N(t)\big), t\in[0,T]\}$ is a Markov process only if this pair is considered as an entity. In other words, a single component, either $\{X_{[q^*, \tilde q]}^N(t), t\in[0,T]\}$ or $\{\tilde M^{tag, \tilde q}_N(t), t\in[0,T]\}$, is not a Markov process and cannot be discussed separately, since the distribution evolution of the $N$ interacting inspectees $\{X_{[q^*, \tilde q]}^N(t), t\in[0,T]\}$ is coupled with any individual's dynamics $\{\tilde M^{tag, \tilde q}_N(t), t\in[0,T]\}$.

In the continuum limit, since any single inspectee's behaviour has negligible impact on the whole population's statistical behaviour,
the distribution dynamics $\{X_{[q^*, \tilde q]}(t), t\in[0,T]\}$ is still the solution to the ordinary differential equation \eqref{limit distribution}, although one inspectee chooses a different strategy from $q^*$. Since $\{X_{[q^*, \tilde q]}(t), t\in[0,T]\}$ is a deterministic process and not affected by any single inspectee's behavior $\{\tilde M^{tag, \tilde q}(t), t\in[0,T]\}$, one can view $\{\tilde M^{tag, \tilde q}(t), t\in[0,T]\}$ as a Markov process, parameterised by $\{X_{[q^*, \tilde q]}(t), t\in[0,T]\}$.
\end{remark}

We show that, as $N\to \infty$, the Markov process $\{\big(X_{[q^*, \tilde q]}^N(t), \tilde M^{tag, \tilde q}_N(t)\big), t\in[0,T]\}$ generated by $\widehat L_t ^N$ in \eqref{tagged L} converges to the Markov process \newline
$\{\big(X_{[q^*, \tilde q]}(t), \tilde M^{tag, \tilde q}(t)\big), t\in[0,T]\}$ generated by $\widehat A_t$ in \eqref{tagged A}. This result  is crucial for  the final result which is stated in Theorem \ref{e convergence}.

To this end, we need the concept of propagators. For a set of continuous function $C(\Sigma_d, \mathbb{L}_d)$, a family of mappings $\Psi^{t,r}$ from $C(\Sigma_d, \mathbb{L}_d)$ to itself, parametrized by the pairs of numbers $r\leq t$ (resp. $t\leq r$) from a given finite or infinite interval is called a {\it (forward) propagator} in $S$, if $\Psi^{t,t}$ is the identity operator in $C(\Sigma_d, \mathbb{L}_d)$ for all $t$ and the following {\it chain rule}, or {\it propagator equation}, holds for $r\leq s\leq t$:
$$\Psi^{t,s}\Psi^{s,r}=\Psi^{t,r}.$$

Let $ \Psi_{N; tag}^{0,t}[q^*, \tilde q]$ denote the propagator generated by  $\widehat L_t ^N[q^*, \tilde q]$ in \eqref{tagged L} and $\Phi_{tag}^{0,t}[q^*, \tilde q]$ the propagator generated by  $\widehat A_t [q^*, \tilde q]$ in  \eqref{tagged A}. By saying this, we mean that for $ f\in \mathbb{D} (\widehat L_t ^N)$, the equations
\begin{equation*}
       \frac{d}{ds}\Psi_{N; tag}^{t,s}f = \Psi_{N; tag}^{t,s}\widehat L_s ^Nf, \quad \frac{d}{ds}\Psi_{N; tag}^{s,r}f = -\widehat L_s ^N\,\Psi_{N; tag}^{s,r}f, \quad 0\leq t\leq s\leq r,
\end{equation*}
hold a.s. in $s$ and for $ f\in \mathbb{D} (\widehat A_t)$, the equations
\begin{equation*}
       \frac{d}{ds}\Phi_{tag}^{t,s}f = \Phi_{tag}^{t,s}\,\widehat A_sf, \quad \frac{d}{ds} \Phi_{tag}^{s,r}f = -\widehat A_s \,\Phi_{tag}^{s,r}f, \quad 0\leq t\leq s\leq r,
\end{equation*}
hold a.s. in $s$.

\begin{theorem}
\label{convergence of propagators}

Suppose that   as $N\to \infty$, the initial data $x_0^N\in \Sigma_d$ converges to certain  $x_0\in \Sigma_d$ and the initial data $m_N(0)\in \mathbb{L}_d$ converges to certain  $m(0)\in \mathbb{L}_d$. Then  for any switching function $\tilde q$ which is Lipschitz continuous in $x$  and any $f\in [C^1(\Sigma_d\times\mathbb{L}_d)]^d$
\begin{equation}
\label{converge}
\lim_{N\to \infty}\left|\Psi_{N; tag}^{0,t}[q^*, \tilde q]  f(x_0^N, m_N(0)) - \Phi_{tag}^{0,t}[q^*, \tilde q]  f(x_0,m(0))\right| =0
\end{equation}
uniformly in $t\in[0,T]$ with any $T>0$.
\end{theorem}

\proof

First, if the space $[C^1(\Sigma_d\times\mathbb{L}_d)]^d$ is a core of the propagator $\Phi_{tag}^{0,t}[q^*, \tilde q]$  generated by the limiting operator  $\widehat A_t [q^*, \tilde q]$ defined in \eqref{tagged A}, then  we have the convergence of the generators on the core of the liming semigroup, i.e., for any $f\in[C^1(\Sigma_d\times\mathbb{L}_d)]^d$
$$\lim_{N\to\infty} \widehat L_t ^N[q^*, \tilde q]f(x, l_k)=\widehat A_t [q^*, \tilde q]f(x, l_k)$$
which implies the convergence of the semigroup on the core $[C^1(\Sigma_d\times\mathbb{L}_d)]^d$, c.f. Kallenberg (2002), i.e. the statement \eqref{converge} is proved.

The result that $[C^1(\Sigma_d\times\mathbb{L}_d)]^d$ is a core of the propagator $\Phi_{tag}^{0,t}[q^*, \tilde q]$  generated by the limiting operator  $\widehat A_t [q^*, \tilde q]$ defined in \eqref{tagged A} is proved in Appendix \ref{Appendices} in 3 steps, see Appendix \ref{A1}-Appendix \ref{A3}.

Here we only need to check the conditions in Appendix \ref{A3} are satisfied. By \eqref{optimal q 1} and  \eqref{u-Lip} we have that the optimal switching function $q^*(t,l_i,x)$ is Lipschitz continuous in $x$. It is clear that
from \eqref{limit distribution}, the function $F$ in \eqref{A3 switching q} with $F_i(t,x)= \sum _{j=1}^d x_j q^*_{ji}(t,x;  V(t,\cdot) )$ is Lipschitz continuous in both $t$ and $x$. Together with the condition that $\tilde q$ is Lipschitz in the variable $x$, the conditions in Appendix \ref{A3} are satisfied. The proof is completed.\qed

\begin{remark}
This kind of convergence results of $N$-particle approximations have been proved  e.g. by Kolokoltsov, Troeva and Yang  (2014)  for a diffusion process and by Kolokoltsov, Li and Yang (2012)  for rather general Markov processes with smoothing property (excluding the present case).
\end{remark}

As a direct consequence of Theorem \ref{Thm1} and Theorem \ref{convergence of propagators} , we have the result in the following theorem, stating that any solution derived from the limiting model \eqref{limit HJB}-\eqref{limit distribution} can be used to approximate the one for an $N$ player game.

 \begin{theorem}
 \label{e convergence}
Suppose that

(i) as $N\to \infty$, the initial data $x_0^N\in \Sigma_d$ converges to certain  $x_0\in \Sigma_d$ and the initial data $m_N(0)\in \mathbb{L}_d$ converges to certain  $m(0)\in \mathbb{L}_d$.

(iii) the terminal cost function $J_T: \mathbb{L}_d\times\Sigma_d\to \R$ is Lipschitz in both variables.

Then a  strategy profile
$$\left\{\alpha^*(X(t)), q^*(t, \cdot, X(t); V(t,\cdot)),\dots, q^*(t,\cdot, X(t); V(t,\cdot))\right\}$$
with $\alpha^*$ and $q^*$ defined respectively in \eqref{limiting-Inspector-payoff-1}  and \eqref {optimal q 1}, and with $(X,V)$ being
a solution to the system \eqref{limit HJB}-\eqref{limit distribution}, is an $\epsilon$-equilibrium in any $N$ inspectee inspection game with $\epsilon=\epsilon(N)\to 0$ as $N\to \infty$.
 \end{theorem}

\proof

First denote by  $M_N^{q^*}(t)$ and $M^{q^*}(t)$ the state dynamics of an inspectee with the switching strategy $q^*$ in the finite $N$ inspectees setting and in the continuum liming setting, respectively. Similarly, denote by $X^N_{[q^*]}(t)$ and $X_{[q^*]}(t)$ the state dynamics of the population with every inspectee applying the switching strategy $q^*$ in the finite $N$ inspectees setting and in the continuum liming setting, respectively. In fact, $X_{[q^*]}(t)$ and $X_{[q^*, \tilde q]}(t)$ generated by the operator $\widehat A_t[q^*, \tilde q]$ in \eqref{tagged A} are the same object, namely the solution to Eq. \eqref{limit distribution}, see Remark \ref{Markov processes}.

To show that the strategy $\alpha^*(X_{[q^*]}(t))$ with $X_{[q^*]}(t)$ being a solution to \eqref{limit distribution} is an $\epsilon$-equilibrium for the inspector, we need to show that,  for $\ep=\ep(N) >0$,
\begin{equation} \label{step 2}
 U_N(a^*(t), X_{[q^*]}(t))\geq U_N(a_N^*(t), X_{[q^*]}^N(t))-\epsilon
 \end{equation}
where the payoff function $U_N$ defined in Eq \eqref{Inspector-payoff-1}, and the inspector's best response functions $ \alpha_N^*(t)$ and $\alpha^*(t)$ are defined in Eq  \eqref{Inspector best response} and Eq \eqref{limiting-Inspector-payoff-1} in the $N$ inspectees setting and  in the mean field inspection setting, respectively.

By the definition of $ \alpha_N^*(t)$ and $\alpha^*(t)$, we have that $ \alpha_N^*(t)$ and $\alpha^*(t)$ are Lipschitz continuous in $X^N(t)$ and $X(t)$ respectively.
Then by Theorem \ref{convergence of propagators} we have
$$\lim_{N\to \infty}| \alpha^*(t)- \alpha_N^*(t)|=0. $$
That is, for $N$ big enough, we have
$$  \alpha^*(t)= \alpha_N^*(t)\pm\epsilon$$
with $\epsilon\to 0$ as $N\to\infty$.
Therefore for $N$ big enough, again by Theorem \ref{convergence of propagators}
\begin{align*}
&U_N( \  \alpha^*(t), X_{[q^*]}(t) )- U_N(   \alpha_N^*(t), X_{[q^*]}^N(t) )\\
= & \E \left (-\alpha^*(t)+L \sum_{i=1}^d   X_{i,[q^*]}(t)\Big(P(\alpha^*(t)) \sigma l_i - (1-P(\alpha^*(t))) l_i\Big)\right)\\
&- \E \left (-\alpha_N^*(t)+L \sum_{i=1}^d   X_{i,[q^*]}^N(t)\Big(P(\alpha_N^*(t)) \sigma l_i - (1-P(\alpha_N^*(t))) l_i\Big)\right)\\
>&-\epsilon,
\end{align*}
with $\epsilon\to 0$ as $N\to\infty$, where the continuity of $P$ in $\alpha$ is used.

Next, to show that the strategy $q^*$ is an $\epsilon$-equilibrium for an individual inspectee, we need to show,  for any inspectee $a=1, \dots, N$ and  $\ep =\ep(N)>0$
\begin{equation} \label{step 1}
J^{(a)}(t,l_i, q^*; X_{[q^*]}^N)\geq J^{(a)}(t,l_i, \tilde q; X_{[q^*,\tilde q]}^N)-\epsilon
\end{equation}
for any $\tilde q$ where the payoff function $J^{(a)}$ is defined in \eqref{insepctee payoff}.
Since the payoff function $J^{(a)}$ defined in \eqref{insepctee payoff} is Lipschitz in $X^N$,  the appendix \ref{A4} implies Eq. \eqref{step 1}.
The proof is completed.
\qed

\section{Approximation analysis}

In this section, we will approximate the optimal switching strategy $q^*$ defined in \eqref{optimal q 1} by a sequence of smooth function $q^*_{\eta}$, for $\eta >0$. We prove that  any smooth approximation $q^*_{\eta}$ as a solution to a forward-backward model is also $\epsilon$ -Nash to any inspection game with finite-number inspectees. This approximation analysis is motived by the following two considerations.

Firstly, based on the result that there exists a solution, a consistent pair $(X(t), V(t,\cdot))$, to the system of equations  \eqref{limit HJB}-\eqref{limit distribution}, one can obtain an optimal investment strategy $\alpha^*(X(t))$  by \eqref{limiting-Inspector-payoff-1} and an optimal switching strategy $q^*(t, l_i, X(t); V(t,\cdot))$ by \eqref{optimal q 1}.  Since there are no analytic formulae for computing $X$ and $V$, one needs to find numerical solutions of $X$ and $V$. In this paper, we do not attempt to investigate methods for abstaining numerical solutions but we intend to adapt our results proved in the previous sections for numerical analysis.
Recall that the obtained optimal switching strategy $q^*$ defined in \eqref{optimal q 1} is a Lipschitz function in $x$. However, for numerical analysis, very often the smoothness of the function $q^*$ is needed. The results in Theorem \ref{smooth strategy} makes our theoretical results in section 4 applicable to numerical analysis.

Secondly, considering smooth approximations $q^*_{\eta}$ as $\epsilon$ -Nash equilibria to the inspection game with finite-number inspectees enables us to discuss the error bound of the approximation as $N\to\infty$.


Now by standard procedure we construct a sequence of matrix-valued smooth functions $q_\eta^*:[0,T]\times \mathbb{L}_d\times\Sigma_d\to \R^d$ to approximate the continuous function $q^*$. We define for $\eta >0$
\begin{equation}
\label{g eta}
q_\eta^*(t,l_i, x):=\int_{\Sigma_d} q^*(t, l_i,x-y)\phi_\eta(y)dy
\end{equation}
where the function $\phi_\eta$ is a smooth mollifier. We have that for any $t\in[0,T]$ and $l_i\in\mathbb{L}_d$, $q^*_{\eta}$ converges to $q^*$ uniformly on $\Sigma_d$, i.e.
$$\lim_{\eta\to 0}\sup_{x\in \Sigma_d} |q^*_\eta(t,l_i,x)- q^*(t,l_i,x)|=0.$$
A typical example of the molllifier $\phi_\eta $ can be $\phi_\eta (y)=\frac{1}{\sqrt{2\pi\eta}}e^{-\frac{y^2}{2\eta}}$.

\begin{theorem}
\label{smooth strategy}
Suppose that as $N\to \infty$, the initial data $x_0^N\in \Sigma_d$ converges to certain  $x_0\in \Sigma_d$ and  the initial data $m_N(0)\in \mathbb{L}_d$ converges to certain  $m(0)\in \mathbb{L}_d$. Moreover,   the terminal cost function $J_T: \mathbb{L}_d\times\Sigma_d\to \R$ is Lipschitz in both variables. Then

(i) any $q_\eta^*$ defined in \eqref{g eta} is an $\epsilon$-Nash for a finite game with
$\epsilon = \epsilon(\eta, N)\to 0$ as $N\to \infty$ and $\eta\to 0$.

(ii) if $q_\eta^*$ defined in \eqref{g eta} is two continuously differentiable in $x$, $\epsilon$ is of order $1/N$.

\end{theorem}
\proof
 (i) To prove $q_\eta^*$ is an $\epsilon$-Nash, we aim to prove that for $(t,l_i)\in[0,T]\times\mathbb{L}_d$ and for any other $\tilde q$
 \begin{equation}
 \label{approx q}
 J^{(a)}(t,l_i, q^*_\eta, X_{q^*_\eta}^N)>J^{(a)}(t,l_i,\tilde q, X_{[q^*_\eta, \tilde q]}^N)- \epsilon.
 \end{equation}
Take the approximating smooth optimal control $q_\eta^*(t,l_i,x)=(q^*_{\eta, i1}(t,X_\eta(t)), \dots, q^*_{\eta, id}(t,X_\eta(t)))$ and consider the system
 \begin{align}
 \label{Xeta}
\left\{
\begin{array}{c l}
&\frac{dX_{\eta,i}(t)}{dt}=\sum _{j=1}^d X_{\eta,j}(t) q^*_{\eta, ji}(t,X_\eta(t)), \quad i=1 ,\dots, d \\[.6em]
& X(0)=x(0)
\end{array}
\right.
\end{align}
Let $X_{q^*_\eta}$ be the solution to the system \eqref{Xeta}. We have that $X_{q^*_\eta}$ converges to $X_{q^*}$ as $\eta \to 0$. By Theorem \ref{convergence of propagators}, as $N\to \infty$, $X^N_{q^*_\eta}\to X_{q^*_\eta}$ and $X^N_{q^*}\to X_{q^*}$. Therefore we have $X^N_{q^*_\eta}\to X^N_{q^*}$ as $N\to \infty$ and $\eta\to 0$.
 Since the payoff function $J^{(a)}$ defined in \eqref{insepctee payoff} is continuous and Lipschitz continuous in $X^N$, we have  for small enough $\eta>0$ and big enough $N$
\begin{align*}
J^{(a)}(t,l_i,q^*_\eta, X_{q^*_\eta}^N)&=J^{(a)}(t,l_i,q^*_\eta, X^N_{q^*})+J^{(a)}(t,l_i,q^*_\eta, X_{q^*_\eta}^N)-J^{(a)}(t,l_i,q^*_\eta, X_{q^*}^N)\\
&=J^{(a)}(t,l_i,q^*_\eta, X^N_{q^*})\pm \epsilon(\eta, N, q^*)\\
&=J^{(a)}(t,l_i,q^*, X^N_{q^*})\pm \epsilon(\eta, N, q^*)\\
&>J^{(a)}(t,l_i,q^*, X^N_{q^*})-\epsilon(\eta, N, q^*)\\
&\geq J^{(a)}(t,l_i,\tilde q, X^N_{[q^*_{\eta}, \tilde q]})-\epsilon(\eta, N, \tilde q, q^*)
\end{align*}
 where the result in the appendix \ref{A4} is used. Hence \eqref{approx q} is proved.

 (ii) To prove $\epsilon$ is of order $1/N$ for a twice continuously differentiable $q^*_\eta$, we aim to show that for $f\in [C^2(\Sigma_d\times\mathbb{L}_d)]^d$,
 \begin{align}
 \label{1/N}
&\left|\Psi_{N; tag}^{0,t}[q^*_\eta, \tilde q]  f(x_0^N, m_N(0)) - \Phi_{tag}^{0,t}[q^*_\eta, \tilde q]  f(x_0,m(0))\right| \leq C(T)\frac{1}{N}\|f\|_{[C^2(\Sigma_d\times\mathbb{L}_d)]^d}.
\end{align}
Here $ \Psi_{N; tag}^{0,t}[q^*_\eta, \tilde q]$ denote the propagator generated by  $\widehat L_t ^N[q^*_\eta, \tilde q]$ and $\Phi_{tag}^{0,t}[q^*_\eta, \tilde q]$ the propagator generated by  $\widehat A_t [q^*_\eta, \tilde q]$. The space $[C^2(\Sigma_d\times \mathbb{L}_d)]^d$ is the set of continuous and bounded vector-valued functions $f$ on $\Sigma_d\times\mathbb{L}_d$ which are twice continuously differentiable in the first variable.

By \eqref{tagged L} and \eqref{tagged A}, we write
\begin{align}
 \label{tagged L eta}
\widehat L_t ^N[q^*_\eta, \tilde q] f(y,l_k)=&\sum_{\substack{i,j=1 \\ i\neq k}}^d (y\cdot e_i) q^*_{\eta,ij}(t,y) N \left [   f(y-\tfrac{e_i}{N}+\tfrac{e_j}{N}, l_k)-f(y,l_k)\right]\notag\\
&+ \sum_{j=1}^d (y\cdot e_k-\tfrac{1}{N}) q^*_{\eta,kj}(t,y) N \left [   f(y-\tfrac{e_{k}}{N}+\tfrac{e_j}{N}, l_k)-f(y,l_k)\right]\notag\\
&+\sum_{j=1}^d \tfrac{1}{N} \tilde q_{kj}(t,y) N \left [   f(y-\tfrac{e_k}{N}+\tfrac{e_j}{N}, l_j)-f(y,l_k)\right]
\end{align}
and
\begin{align}
\label{tagged A eta}
\widehat A_t[q^*_\eta, \tilde q]f(x,l_k)
=& \sum _{\substack{j=1 \\ i\neq j}}^d\left(x_iq^*_{\eta, ij}(t,x) -x_j q^*_{\eta, ji}(t,x)\right)  \frac{\partial f}{\partial x_j }(x,l_k)\notag\\
&+ \sum _{\substack{j=1 \\ i\neq j}}^d \tilde q_{kj(t,x)}(f(x,l_j)-f(x,l_k)).
\end{align}
In fact, since $q^*_\eta$ is twice continuously differentiable in $x$, the space $[C^2(\Sigma_d\times \mathbb{L}_d)]^d$ is a invariant core for $\widehat L_t ^N[q^*_\eta, \tilde q]$ and $\widehat A_t[q^*_\eta, \tilde q]$,  that is
$$\widehat L_t ^N[q^*_\eta, \tilde q]: [C^2(\Sigma_d\times \mathbb{L}_d)]^d\to [C^2(\Sigma_d\times \mathbb{L}_d)]^d$$
and
$$\widehat A_t [q^*_\eta, \tilde q]: [C^2(\Sigma_d\times \mathbb{L}_d)]^d\to [C^2(\Sigma_d\times \mathbb{L}_d)]^d.$$
Further, by Taylor theorem, \eqref{tagged L eta} can be expended by using the following representation (ref. Kolokoltsov (2010), Corollary F.2)
\begin{align}
 &f(y-\zeta, l_k)-f(y,l_k)=\left(\frac{\partial f(y,l_k)}{\partial \zeta},\zeta\right) +\int_0^1 ds (1-s)\left( \frac{\partial^2 f(y,l_k)}{\partial \zeta^2},\zeta^2 \right).
\end{align}
 Therefore, for $f\in[C^2(\Sigma_d\times \mathbb{L}_d)]^d$, we have
 \begin{align}
 \label{A N}
 \|(\widehat L_t ^N[q^*_\eta, \tilde q]- \widehat A_t[q^*_\eta, \tilde q])f\|_{[C^2(\Sigma_d\times \mathbb{L}_d)]^d}\leq C(T) \frac{1}{N}\|f\|_{[C^2(\Sigma_d\times\mathbb{L}_d)]^d}.
  \end{align}
 To complete the proof, we need the following calculation: for $s\leq t$
 \begin{align}
 \label{trick}
 \Psi_{N; tag}^{s,t}  f - \Phi_{tag}^{s,t} f
& =\Psi_{N; tag}^{s,r}\Phi_{tag}^{r,t}[|_{r=s}^tf=\int_s^t\frac{d}{dr}\left(\Psi_{N; tag}^{s,r} \Phi_{tag}^{r,t}\right) f dr\notag\\
 &= \int _s^t \Psi_{N; tag}^{s,r}  \left(  \widehat L_t ^N- \widehat A_t    \right)\Phi_{tag}^{r,t}f dr.
 \end{align}
By \eqref{A N} and \eqref{trick} together, we get the required statement \eqref{1/N}. Consequently, by the definition of $J^{(a)}$ in \eqref{insepctee payoff} we have
\begin{align*}
\sup_{t}\left|   J^{(a)}(t,\cdot,q^*_\eta, X_{q^*_\eta}^N)      -  J^{(a)}(t,\cdot,\tilde q, X^N_{[q^*_{\eta}, \tilde q]})\right|\leq C(T)\frac{1}{N}\|f\|_{[C^2(\Sigma_d\times\mathbb{L}_d)]^d}.
\end{align*}
\qed

\section{Appendix}
\appendix
\section{Proof to Theorem \ref{Thm1}}
\label{Proof to Theorem Thm1}
The proof to Theorem \ref{Thm1} consists of three steps.

Step 1: for any given $ X\in C([0,T], \Sigma_d)$, we show that the HJB equation \eqref{limit HJB} is well posed. Moreover, the resulting solution, denoted by $V(t,l_i;X)$, is Lipschitz with respect to the parameter $X$.

For proving the existence of a solution to the ordinary differential equation \eqref{limit HJB}, it is sufficient to have that
the function $H$ defined in \eqref{H-def}
is Lipschitz in $\phi$ uniformly.

Since the optimal switching function $q^*_{ij}$ in \eqref{optimal q 1} is Lipschitz continuous in $\phi$, we conclude that there exists a constant $c>0$ such that for any $l_i\in\mathbb{L}_d, x\in\Sigma_d$,
\begin{equation}
\label{H lip in phi}
|H(l_i,\phi,x)- H(l_i,\psi,x)|\leq c \|\phi-\psi\|, \quad \forall  \phi,\psi\in\R^d.
\end{equation}
Therefore, for any $ X\in C([0,T], \Sigma_d)$,  there exists a unique solution to \eqref{limit HJB}.

To show that the solution $V$ is Lipshitz with respect to the parameter $X$, we write the equation \eqref{HJB-g} in integral form:
$$V(t,l_i;X)+\int_t^T H(l_i, V(s,\cdot;X),X(s))ds=0.$$
The function $H$ defined in \eqref{H-def}
is Lipschitz in $x$ uniformly, since for any $ l_i\in\mathbb{L}_d, \phi\in\R^d$ and any $x, y \in \Sigma_d$,
\begin{align}
\label{H lip in mu}
|H(l_i,\phi,x)- H(l_i,\phi,y)|&\leq l_i(1+\sigma)  \left| P(\alpha^*(x))-P( \alpha^*(y))  \right|\leq c\|x-\eta\|_{TV}
\end{align}
with a constant $c>0$.
For any $t\in[0,T], l_i\in\mathbb{L}_d$ and $X,Y\in C([0,T], \Sigma_d)$
\begin{align}
&|V(t,l_i;X)-V(t,l_i;Y)|\notag\\
&\leq \int _t^T |H(l_i, V(s,\cdot;X),X(s))-H(l_i, V(s,\cdot;Y),Y(s)) |ds\notag\\
&\leq \int _t^T |H(l_i, V(s,\cdot;X),X(s))-H(i, V(s,\cdot;X),Y(s)) |ds \notag\\
&\hspace{1.5em}+\int _t^T |H(l_i, V(s,\cdot;X),Y(s))-H(l_i, V(s,\cdot;Y),Y(s)) |ds\notag\\
&\leq c T \|X-Y\|_{\infty} + c \int_t^T \|V(s,\cdot;X)-V(s,\cdot;Y)\| ds\notag
\end{align}
where \eqref{H lip in phi} and \eqref{H lip in mu} are used to get the last inequality. Then by Gronwall's inequality,  the solution to \eqref{limit HJB} is Lipschitz continuous in $X$, i.e. there exists a constant $c>0$ such that for any $t\in[0,T]$,
\begin{align*}
\|V(t,\cdot;X)-V(t,\cdot;Y)\|\leq c T\|X-Y\|_{\infty}.
\end{align*}
Consequenctly,  by the definition of $q^*$ in \eqref{optimal q 1}, we have that for any $t\in[0,T]$ and $i,j=1,\dots,d$ with $i\neq j$,
\begin{equation}
\label{u-Lip}
|q^*_{ij} (t,X(t); V(t,\cdot;X))-q^*_{ij}(t,Y(t); V(t,\cdot;Y)) |\leq cT\|X-Y\|_{\infty}.
\end{equation}

Step 2: we prove the existence of the solution to the following equation
\begin{align}
 \label{limit distribution 2}
\left\{
\begin{array}{c l}
&\frac{dX_i(t)}{dt}=\sum _{j=1}^d X_j(t) q_{ji}(t,X(t)), \quad i=1 ,\dots, d \\[.6em]
& X(0)=x(0)
\end{array}
\right.
\end{align}
with a given switching policy $q$ which is Lipschitz in $x$ and then prove the sensitivity of the solution $X(\cdot)$ with respect to those $q$ which are Lipschitz in $x$.

 Define a vector field $G: [0,T] \times\Sigma_d\times\R\to \R^d$ with its $i$th component $G_i:  [0,T]\times\Sigma_d\times\R\to \R$
 \begin{align}
 G_i(t,x,q):=\sum _{j=1}^d x_j q_{ji}(t, x).
 \end{align}
To prove the existence of a solution to \eqref{limit distribution} it is sufficient to prove that  $G$ is Lipschitz continuous in $x$. We say $G$ is Lipschitz continuous in $x$ if each component $G_i$, $i=1,\dots, d$, is Lipschitz in $x$.
 Now we show that $G_i$ is Lipschitz in $x$ by using Schwarz inequality. Since $q$ is assumed to be Lipschitz in $x$, we have for  $t\in[0,T]$ and $x, y\in\Sigma_d$,
  \begin{align}
  \label{G-Lip}
  | G_i(t,x,q)- G_i(t,y,q)|&\leq |\sum _{j=1}^d x_j  (q_{ji}(t, x) -q_{ji}(t, y) ) | + | \sum _{j=1}^d(x_j-\eta_j)  q_{ji}(t, y)    |\notag\\
  &\leq c \|x-y\|+ Q \|x-y\|
     \end{align}
 which implies that  $G_i$ is Lipschitz in the second variable $x$. Hence, for any given $q$ which is Lipschitz in $x$, there exists a solution to \eqref{limit distribution 2}.

 To prove the sensitivity of the solution to \eqref{limit distribution 2} with respect to those $q$ which are Lipschitz in $x$, we rewrite \eqref{limit distribution 2}  in a integral form.
Let $X$ ( resp. $Y$) be the solution to \eqref{limit distribution 2}  under the control policy $q_1$ (resp. $q_2$), with the same initial value $x_0\in\Sigma_d$, namely
\begin{align*}
   X_i(t)=x_{0,i}+\int_0^t G_i(s, X(s), q_1(s, l_i, X(s)))ds
    \end{align*}
    and
     \begin{align*}
   Y_i(t)=x_{0,i}+\int_0^t G_i(s, Y(s), q_2(s, l_i, Y(s)))ds.
    \end{align*}
Then by \eqref{G-Lip}, we have
       \begin{align}
   |X_i(t)-Y_i(t)|&\leq \int_0^t |G_i(s, X(s), q_1(s, l_i, X(s)))ds-  G_i(s, Y(s), q_2(s, l_i, Y(s)))| ds\notag\\
   &\leq \int_0^t |G_i(s, X(s), q_1(s, l_i, X(s)))- G(s,X(s), q_2(s, l_i, X(s)))|ds\notag\\
   & \hspace{1em}+ \int_0^t |G_i(s, X(s), q_2(s, l_i, X(s)))-G_i(s, Y(s), q_2(s, l_i, Y(s)))|ds\notag\\
   &\leq \int_0^t c|q_1(s, l_i, X(s))- q_2(s, l_i, X(s))|ds+ \int_0^t c \|X(s)-Y(s)\|_{TV}ds.\notag
    \end{align}
    By Gronwall's inequality, we get for any $t\in[0,T]$
    \begin{align}
    \label{X-Lip}
    \|X(t)-Y(t)\|_{TV}&=\sum_{i=1}^d  |X_i(t)-Y_i(t)|   \leq c T\sup_{l_i\in\mathbb{L}_d} \|q_1(\cdot, l_i, X(s))-q_2(\cdot, l_i, X(s))\|_{\infty}
    \end{align}
    with a constant $c>0$.

Step 3: In summary, so far we have considered the following  mapping
\begin{align}
\label{mapping}
&\quad X\quad\quad\to\quad \quad q^* \quad \to\quad \quad\bar X\notag\\[.4em]
\Gamma: \quad &C([0,T], \Sigma_d)\quad \, \,\Longrightarrow \quad \quad C([0,T], \Sigma_d).
\end{align}
By analysing \eqref{limit HJB}  we get \eqref {u-Lip}, namely the resulting  optimal switching function $q^*(t,l_i,x) $  is Lipschitz with respect to $x$; further, by analysing \eqref{limit distribution} with any switching policy $q$ which is Lipschitz in $x$, we get \eqref{X-Lip}, namely the solution to  \eqref{limit distribution}  is Lipschitz with respect to its control parameter $q$.
Therefore we can conclude that the mapping $\Gamma$  \eqref{mapping} is Lipschitz, that is, for any $X,Y\in C([0,T], \Sigma_d)$,  there exists a constant $c>0$ such that
\begin{equation}
\|\Gamma (X)-\Gamma (Y)\|_{\infty}=\|\bar X-\bar Y\|_{\infty}\leq c T \|  X-Y\|_{\infty}
\end{equation}
where the norm $\|\cdot\|_{\infty}$ is defined in \eqref{norm}.

Thus for a small $T$, the mapping $\Gamma$ is a contraction, proving statement (i). For arbitrary finite $T$,  one has, from \eqref{Kinetic}, that the image of the mapping $\Gamma$ is bounded equicontinuous and hence a compact subset of $C([0,T], \Sigma_d)$ (by Arzela-Ascoli  theorem). Hence by the Brouwer fixed point theorem, $\Gamma$ has a fixed point, proving statement (ii).

\section{$[C^1(\Sigma_d\times\mathbb{L}_d)]^d$  is a core for the generator
$\widehat A_t$}
\label{Appendices}
This appendix aims to prove step by step that the space $[C^1(\Sigma_d\times\mathbb{L}_d)]^d$ is a core for the limiting generator
$\widehat A_t$ defined in \eqref{tagged A}.  In \ref{A1}, we consider a single deterministic system $X(t), t\in[0,T]$ and show that $[C^1(\Sigma_d)]^d$ is a core for the generator of the system. Then in \ref{A2}, we consider a time-homogenous Markov chain $M(t), t\in[0,T]$ modulated by a deterministic system $X(t), t\in[0,T]$. We show that $[C^1(\Sigma_d\times\mathbb{L}_d)]^d$ is a core for the generator of the system $(M(t), X(t)), t\in[0,T]$. Finally in \ref{A3}, we consider a time- nonhomogenous Markov chain $M(t), t\in[0,T]$ modulated by a deterministic system $X(t), t\in[0,T]$. This Fellow process $(M(t), X(t)), t\in[0,T]$ is exactly the one generated by the limiting operator $\widehat A_t$ defined in \eqref{tagged A}. We show that $[C^1(\Sigma_d\times\mathbb{L}_d)]^d$ is indeed a core for the generator of the system $(M(t), X(t)), t\in[0,T]$, namely the generator
$\widehat A_t$ defined in \eqref{tagged A}.

Recall that the space $[C^1(\Sigma_d\times \mathbb{L}_d)]^d$ is the set of continuous and bounded vector-valued functions $f$ on $\Sigma_d\times\mathbb{L}_d$ which are differentiable in the first variable.
%
%
The standard notation $\dot {}$ denotes the differentiation with respect to time, e.g. $\dot x=\frac{dx}{dt}$.
\subsection{Evolution of the deterministic dynamics $X(t)$}
\label{A1}
Consider a system $X(t), t\in[0,T]$ which is described by the first-order ordinary differential equation
\begin{equation} \label{case 1}
\dot X_t=F(X_t)
\end{equation}
with a given initial $x_0\in \Sigma_d$ and $F:\Sigma_d\to \R^d$. The system  $X(t), t\in[0,T]$  has the generator $A:[C^1(\Sigma_d)]^d\to [C(\Sigma_d)]^d$ which is of the form
 \begin{align}
\label{limit operator A1}
Af(x) = F(x) \frac{\partial f}{\partial x}(x).
\end{align}
Let $\phi_t$ denote the semigroup generated by the generator $A$ in \eqref{limit operator A1}. The solution of Eq. \eqref{case 1} is
$$f(X_t(x_0))=(\phi_t f)(x_0).$$
\begin{lemma}
If the function $F$ in \eqref{case 1} is Lipschitz, then $[C^1(\Sigma_d)]^d$ is a core of the generator $A$ in \eqref{limit operator A1}.
\end{lemma}
\proof
Since the space $[C^1(\Sigma_d)]^d$ is not invariant under the operator $A$ \eqref{limit operator A1}, we cannot apply the standard result (c.f. Kallenberg (2002))  that a dense invariant subset of the domain is always a core. We introduce a subspace by collecting all shifted functions from $[C^1(\Sigma_d)]^d$ and their linear combinations and denote this space by $[\tilde C(\Sigma_d)]^d$
$$\tilde C(\Sigma_d):=\{f_{i}(X_t(x_0))\big| \forall t\in[0,T], f_{i} \in C^1(\Sigma_d) \}$$
 with $ i\in\{1,\dots, d\}$. By its construction, this space $[\tilde C(\Sigma_d)]^d$ is an invariant core for the semigroup $\phi_t$.
In order to prove that $[C^1(\Sigma_d)]^d$ is still a core for the semigroup $\phi_t$, we construct a sequence of  vector-valued smooth functions $F^n$ and for each $n\in\mathbb{N}$, the entry $F^n_{i} \in C^1(\Sigma_d),  i=1,\dots, d$  is defined as
$$ F^n_{i}(x):=\int_{\Sigma_d} F_{i}(x-y)\phi_n(y)dy $$
where the mollifier $\phi_n$ is at least first order differentiable and with compact support. So we have, as $n\to \infty$, $F^n_{i} \in C^1(\Sigma_d)$ converges to $F_{ij}\in \tilde C(\Sigma_d)$ as $n\to \infty$. Then we consider the following approximating systems
\begin{equation*}
\dot X_t=F^n(X_t)
\end{equation*}
with  an initial $x_0\in\Sigma_d$. Let  $X^n_t(x_0)$ denote the solution to this approximating systems. Since $F^n$ converges $F$ as $n\to \infty$, we have $X_t^n(x_0)$ converges to $X_t$ as $n\to \infty$.
Therefore the closure of the subspace $[C^1(\Sigma_d)]^d$ is $[\tilde C(\Sigma_d)]^d$ which is a core of $A$, hence $[C^1(\Sigma_d)]^d$ is a core of $A$.

\qed

\subsection{A homogenous Markov chain $M(t)$ modulated by a deterministic evolution $X(t)$}
\label{A2}
Consider a Markov chain $M(t), t\in[0,T]$ on $\mathbb{L}_d$ with a switching function $q$ on $\Sigma_d$
\begin{gather}
\label{A2 switching q}
x\,\to\,\mathbb{Q}(x)=
\left( \begin{array}{ccc}
q_{11}(x), & \dots, & q_{1d}(x) \\
 & \dots &  \\
 q_{d1}(x), & \dots, & q_{id}(x) \\
 & \dots &  \\
q_{d1}(x), & \dots, & q_{dd}(x)
 \end{array} \right).
 \end{gather}
Let the Markov chain $M(t), t\in[0,T]$ be modulated by a deterministic evolution $X(t), t\in[0,t]$ which is described by the first order  ordinary differential equation
\begin{equation}
\label{A2 X}
\dot X_t=F(X_t)
\end{equation}
with a given initial $x_0\in \Sigma_d$ and $F:\Sigma_d\to \R^d$.
Then the modulated Markov chain is a Feller process on $[C(\Sigma_d\times\mathbb{L}_d)]^d$, which is denoted by $(X(t), M(t)), t\in[0,T]$  and is described by
\begin{align}
\label{case 2}
&\dot h=\mathbb{Q}(X_t(x_0))h
\end{align}
with a given initial function $h_0\in C(\Sigma^d\times\mathbb{L}_d)$.  The solution of \eqref{case 2} is
\begin{equation}
\label{h in 2}
h_t(h_0, x_0)= (\psi_th_0)(x_0)
\end{equation}
where $\psi_t$ is the semigroup of the Feller process $(X(t), M(t))$.
\begin{lemma}
If the functions $F$ in \eqref{A2 X} and $q$ in \eqref{A2 switching q} are Lipschitz, then $[C^1(\Sigma_d\times\mathbb{L}_d)]^d$ is a core for the generator of the semigroup $\psi_t$ in \eqref{h in 2}.
\end{lemma}
\proof
Following the proof for A.1, we construct a sequence of  matrix-valued smooth functions $\mathbb{Q}^n $ with each entry $q^n_{ij}\in C^1(\Sigma_d)$ defined as
$$ q^n_{ij}(x):=\int_{\Sigma_d} q_{ij}(x-y)\phi_n(y)dy,\quad i,j=1,\dots, d   $$
where the mollifier $\phi_n$ is at least first order differentiable and with compact support, so that $\mathbb{Q}^n$ converges to $\mathbb{Q}$ as $n\to \infty$. The solutions to the approximating systems
\begin{equation*}
\dot h=\mathbb{Q}^n(X^n_t(x_0))h
\end{equation*}
are denoted by $h^n_t(x_0, h_0)$, where the sequence $X^n_t(x_0)$ is constructed in the proof for A.1. Since  $X_t^n(x_0)\to X_t(x_0)$ as $n\to \infty$ and $\mathbb{Q}$ is a Lipschitz function, we have $\mathbb{Q}^n(X_t^n(x_0))\to \mathbb{Q}(X_t(x_0))$ as $n\to \infty$. Hence, $h^n_t(x_0, h_0) \to h_t(x_0, h_0)$ as $n\to \infty$.
Therefore the closure of the space $[C^1(\Sigma_d\times\mathbb{L}_d)]^d$ is $[\tilde C(\Sigma_d\times\mathbb{L}_d)]^d$, hence $[C^1(\Sigma_d\times\mathbb{L}_d)]^d$ is a core the generator of the semigroup $\psi_t$ in \eqref{h in 2}. \qed

\subsection{A non-homogenous Markov chain $M(t)$ modulated by a deterministic evolution $X(t)$}
\label{A3}

Consider a time non-homogenous Markov chain $M(t), t\in[0,T]$ on $\mathbb{L}_d$ with a switching function $\mathbb{Q}$ on $[0,T]\times\Sigma_d$
\begin{gather}
\label{A3 switching q}
(t,x)\,\to\,\mathbb{Q}(t,x)=
\left( \begin{array}{ccc}
q_{11}(t,x), & \dots, & q_{1d}(t,x) \\
 & \dots &  \\
 q_{i1}(t,x), & \dots, & q_{id}(t,x) \\
 & \dots &  \\
q_{d1}(t,x), & \dots, & q_{d}(t,x)
 \end{array} \right).
 \end{gather}
Let the time non-homogeneous Markov chain $M(t), t\in[0,T]$ be modulated by a deterministic evolution $X(t)$ which is described by
\begin{equation}\label{time non-homo case}
\dot X_t=F(t, X_t)
\end{equation}
with a given initial $x_0\in \R^d$ and $F: [0,T]\times\Sigma_d\to \R^d$.
Then the modulated non-homogenous Markov chain is a Feller process on $[C(\Sigma_d\times\mathbb{L}_d)]^d$, which is denoted by $(X(t), M(t)), t\in[0,T]$  and is described by
\begin{align}
\label{case 3}
&\dot h=\mathbb{Q}(t, X_t(x_0))h
\end{align}
with a given initial function $h_0\in C(\Sigma^d\times\mathbb{L}_d)$.  The solution of \eqref{case 2} is
\begin{equation}
\label{h in 3}
h_t(h_0, x_0)= (\psi^{t,0}h_0)(x_0)
\end{equation}
where $\psi^{t,s}$, $0\leq s\leq t,$  is the two-parameter semigroup of the Feller process $(X(t), M(t))$.
\begin{lemma}
If the function $F$ in \eqref{time non-homo case} and $q$ in \eqref{A3 switching q} are Lipschitz in both $x$ and $t$, then $[C^1(\Sigma_d\times\mathbb{L}_d)]^d$ is a core for the generator of $\psi^{t,s}$ in \eqref{h in 3}.
\end{lemma}
\proof
Set $y=(x,t)\in \Sigma_d \times [0,T]$. Then Markov chain $M(t)$ is governed by the switching function $q(y)$ and the system \eqref{time non-homo case} is translated to
\begin{align*}
&\dot X_t=F(y)\\
&\dot s=1
\end{align*}
with a initial data $(x_0, 0)$. Then a direct application of the result in the appendix \ref{A2} complete the proof.
\qed

%
%
%
%

\section{$\epsilon$ equilibrium under a general payoff function $J$}
\label{A4}

This appendix states that, for an inspection game with a general payoff function $J$ of inspectees,  any optimal $q^*$ derived from the corresponding mean field inspection game model  is an approximate Nash for any inspectee.

Consider all $N$ inspectees aim to maximise  a general payoff function $J$ as a function of $t,l,q$ and $x$.
Let $X^N=\{X^N(t), t\in[0,T]\}$ be the distribution evolution of the $N$ interacting inspectees among $d$ states in $\mathbb{L}_d$.
Let $X=\{X(t), t\in[0,T]\}$ be a solution to the mean field inspection game and $q^*$ be the resulting optimal switching strategy for a representative inspectee.
If $X^N\to X$ and the payoff function $J$ is Lipschitz uniformly in  $x$, then $q^*$ is an $\epsilon$ equilibrium for an inspection games with any finite $N$ inspectees, namely for any $(t,l_i)\in [0,T]\times \mathbb{L}_d$,
$$J(t,l_i,q^*, X_N)\geq J(t,l_i,\tilde q, X_N)-\epsilon$$
where $\epsilon=\epsilon(N, q^*, \tilde q)\to 0$ as $N\to \infty$.
\proof

Since $X^N\to X$ and $J$ is Lipschitz uniformly in  $x$,  for any $(t,l_i)\in [0,T]\times \mathbb{L}_d $ and an switching function $q$, we have
$$\lim_{N\to\infty} J(t,l_i,q, X_N)= J(t,l_i,q,X).$$

Since inspectees aim to maximise their payoffs and in the limit $N\to\infty$, $q^*$ is the optimal strategy,
for any $(t,l_i)\in [0,T]\times \mathbb{L}_d$ we have
$$J(t,l_i,q^*, X)\geq J(t,l_i,\tilde q, X)$$
for any $\tilde q$. Therefore, for any $(t,l_i)\in [0,T]\times \mathbb{L}_d$ and for $N$ big enough,
there exists an $\epsilon=\epsilon(q^*, \tilde q, N)>0$ so that
\begin{align*}
J(t,l_i,q^*, X_N)&=J(t,l_i,q^*,X)+J(t,l_i,q^*, X_N)-J(t,l_i,q^*, X)\\
&=J(t,l_i,q^*,X) \pm\epsilon (q^*, N)\\
&\geq J(t,l_i,\tilde q,X) \pm\epsilon (q^*, N)\\
&=J(t,l_i,\tilde q, X_N)+J(t,l_i,\tilde q, X)-J(t,l_i,\tilde q, X_N)\pm\epsilon (q^*, N)\\
&= J(t,l_i,\tilde q, X_N )\pm\epsilon (\tilde q, N)\pm\epsilon (q^*, N)\\
&\geq  J(t,l_i,\tilde q, X_N )-\epsilon (q^*, \tilde q, N)
\end{align*}
 with $\epsilon=\epsilon(q^*, \tilde q, N)\to 0$ as $N\to \infty$.
\qed

{\small
\nocite{*}
\bibliographystyle{acm}

\bibliography{biblio}

}

\end{document}